\documentclass[preprint]{elsarticle}

\usepackage[english]{babel}
\usepackage{graphicx}
\usepackage{dcolumn}
\usepackage{bm}

\usepackage{upgreek}
\usepackage{amsmath}
\usepackage{amssymb} 
\usepackage{diagbox}
\usepackage{commath}

\usepackage[utf8]{inputenc}
\usepackage{xcolor}

\usepackage[hidelinks]{hyperref} 

\newcommand{\ii}{\mathrm{i}} 

\renewcommand{\t}{\mathbf}

\usepackage{layouts} 

\begin{document}

\title{Elimination of ringing artifacts by finite-element projection in FFT-based homogenization}

\author[1]{Richard J. Leute}
\address[1]{Department of Microsystems Engineering, University of Freiburg, Georges-K\"ohler-Allee 103, 79110 Freiburg, Germany}

\author[2]{Martin Ladecký}
\address[2]{Faculty of Civil Engineering, Czech Technical University in Prague, Th\'akurova 7, 166 29 Prague 6, Czech Republic}

\author[3]{Ali Falsafi}
\address[3]{Department of Mechanical Engineering, \'Ecole Polytechnique F\'ed\'erale de Lausanne, 1015 Lausanne, Switzerland}

\author[1,4]{Indre Jödicke}
\address[4]{Cluster of Excellence livMatS, Freiburg Center for Interactive Materials and Bioinspired Technologies, University of Freiburg, Georges-K\"ohler-Allee 105, 79110 Freiburg, Germany}

\author[2]{Ivana Pultarová}

\author[2]{Jan Zeman}

\author[3]{Till Junge}

\author[1,4]{Lars Pastewka}
\ead{lars.pastewka@imtek.uni-freiburg.de}

\begin{abstract}
Micromechanical homogenization is often carried out with Fourier-accelerated methods that are prone to ringing artifacts. We here generalize the compatibility projection introduced by Vond\v{r}ejc, Zeman \& Marek [Comput. Math. Appl. 68, 156  (2014)] beyond the Fourier basis. In particular, we formulate the compatibility projection for linear finite elements while maintaining Fourier-acceleration and the fast convergence properties of the original method. We demonstrate that this eliminates ringing artifacts and yields an efficient computational homogenization scheme that is equivalent to canonical finite-element formulations on fully structured grids.
\end{abstract}

\maketitle

\section{Introduction}

Mechanical homogenization seeks the computation of effective mechanical properties, such as homogenized elastic constants or the complete (potentially nonlinear) stress-strain response, given the microstructure of material and constitutive laws for the individual components. Homogenization often operates on periodic ``representative'' volume elements that are characteristic for the material under investigation \cite{Michel1999,Matous2017}. The intrinsic periodicity of such volume elements suggests the use of a spectral Fourier-basis~\cite{boyd_chebyshev_2000} for the efficient numerical solution of the homogenization problem. 

Since the seminal works by Moulinec and Suquet in 1994 \cite{Moulinec1994} and 1998 \cite{Moulinec1998}, (fast) Fourier transform (FFT) based homogenization methods were developed rapidly; see, e.g. Ref.~\cite{Schneider2021} for an authoritative review. The spectral method is here used to solve the partial differential equation (PDE) for a static mechanical equilibrium of a microstructure. The availability of highly optimized FFT implementations (like FFTW~\cite{Frigo2005}) enabled the implementation of efficient and highly parallel spectral solvers which can beat standard implementations of the finite element method (FEM) in speed and accessible system sizes~\cite{Roters2019}.

A major advantage of spectral methods is that it is straightforward to increase the basis set and by this systematically improve the accuracy of the solution. However, because of their global support trigonometric polynomials are not well suited for the solution of PDEs with discontinuous material coefficients or discontinuities in their solution. This leads to a phenomenon known as Gibbs ringing. Gibbs ringing in Fourier spectral methods is well documented, and there are several approaches to suppress the problem~\cite{Gottlieb1997}.

Gibbs ringing persists in FFT-based homogenization methods and was investigated by several authors (see e.g. Refs.~\cite{Mueller1996,Kassbohm2006,Willot2014,Brisard2012, Willot2015,Schneider2016,Khorrami2020,Ma2021}). A typically mitigation strategy is the introduction of discrete derivatives: for example, M\"uller \cite{Mueller1996} used a finite difference discretization or Willot~\cite{Willot2015} and Schneider et al.~\cite{Schneider2016} describe a central-difference scheme on a staggered grid. However, published works have mainly focused on natural extensions of the Moulinec-Suquet scheme that is formulated in terms of a reference problem in the displacement space, i.e. for small-strain elasticity in form of the Navier-Lam\'e equations, that explicitly contain the constitutive law. A solution in Fourier space is then only possible for a homogeneous reference medium. Heterogeneity is superimposed onto this reference solution through Lippman-Schwinger-like approaches for capturing it. Convergence properties of such schemes deteriorate with the strength of the heterogeneity.

In this article we discuss the problem of ringing in FFT-based homogenization methods and present a general solution to it for which we combine the advantages of an FEM-type Galerkin approach with the speed-up of a FFT method. We work in the framework of compatibility projection developed by Vond\-{r}ejc, Zeman, de Geus and coworkers~\cite{Vondrejc2014,Zeman2017,deGeus2017} that considers the deformation gradient as the primary degree of freedom and eliminates the need for an explicit reference medium. Previous work, e.g. Refs.~\cite{ISI:000389105900018,Vidyasagar2017}, has proposed finite differences approaches to reducing ringing artifacts. Different to previous works, we derive a general expression for a discretized operator for the projection on compatible fields and we propose an elimination rather than a mitigation of the problem for certain operator choices. We show that in addition to Gibbs ringing, the solution with projection techniques is prone to ringing artifacts arising from missing degrees of freedom in the formulation of the deformation gradient.

Our suggested improvements to the method are supported by numerical simulations and comparison with analytical results. We suggest a linear finite-element-based projection operator that eliminates all ringing artifacts and discuss the analogy of the modified FFT-based homogenization method to standard FEM. All methods are implemented in the open source code $\mu$Spectre \cite{muspectre2021} and the numerical examples shown in Section~\ref{chap:examples} can easily be reproduced by running the corresponding examples of the software.

\section{Methods}
\label{chap:methods}

\subsection{Compatibility projection}

We investigate microstructures by a representative volume element (RVE) in a periodic simulation cell $\Omega_0$. The microstructure can consist of different phases which are described by arbitrary small- or finite-strain material models. Here and in the following we will denote first-order tensors (vectors) by arrows and second-order tensors by bold symbols.

Any deformation of the simulation cell can be described by a function (the placement map) $\vec{\chi}:  \Omega_0 \rightarrow \Omega$, mapping the undeformed grid positions $\vec{r} \in \Omega_0$ into a deformed configuration $\vec{\chi}(\vec{r}) \in \Omega$ where $\Omega$ is the deformed periodic simulation domain. The overall goal is to solve for the static mechanical equilibrium of the periodic cell for a given deformation. The static mechanic equilibrium is given by \cite{Belytschko2014}
\begin{equation}
    \label{eq:stat-mech-equi}
    \nabla_0\cdot\t{P}^T(\t{F}(\vec{r})) = \vec{0}
    \quad\text{or}\quad
    \partial_{0,\alpha} P_{i\alpha}(\t{F}(\vec{r})) = 0,
\end{equation}
with implicit summation over repeated indices (the Einstein summation convention) and a dot product ($\cdot$) as defined by the second part of Eq.~\eqref{eq:stat-mech-equi}. Greek indices indicate a tensor dimension representing a derivative, whereas Latin indices are used for all other tensor dimensions.
$\t{P}$ is the first Piola-Kirchhoff stress tensor, in general a non-linear function of the deformation gradient
\begin{equation}
    \label{eq:def-grad}
    \t{F}(\vec{r}) = \nabla_0 \vec{\chi}(\vec{r})
    \quad \text{or} \quad
    F_{i\alpha}(\vec{r}) = \partial_{0,\alpha} \chi_{i}(\vec{r})
\end{equation}
Note that the partial derivatives $\partial_{0,\alpha} \equiv \partial/\partial r_\alpha$ in Eqs.~\eqref{eq:stat-mech-equi} and \eqref{eq:def-grad} are with respect to the undeformed configuration $\Omega_0$ of the cell.

The most widely used homogeneization schemes combine Eqs.~\eqref{eq:stat-mech-equi} and \eqref{eq:def-grad}, leading to a set of second-order differential equations in $\vec{\chi}$. For small strains, those are the well known Navier-Lam\'e equations that contain spatial derivatives of the elastic constants. The Green's function of this second order differential equation therefore explicitly contains the constitutive law. Solution with Fourier-techniques then requires the introduction of a homogeneous reference material.

In contrast, the formulation employed here \cite{Vondrejc2014,Zeman2017,deGeus2017} solves the set of first-order differential equations given by Eqs.~\eqref{eq:stat-mech-equi} and \eqref{eq:def-grad}. This leaves the deformation gradient $\t{F}$ as a degree of freedom. Equation~\eqref{eq:def-grad} can be interpreted as the constraint that $\t{F}$ needs to be \emph{compatible}, i.e. given by the gradient of a respective placement map. The mathematical trick is that this constraint can be implemented by a projection operator $\mathbb{G}$ that turns any second-order tensor into a compatible one.

Following Refs.~\cite{Vondrejc2014,Zeman2017,deGeus2017}, we reformulate Eq.~\eqref{eq:stat-mech-equi} in the weak (weighted residual) form,
\begin{equation}
\label{eq:stat-mech-equi-weak-form}
    \int_{\Omega_0} \dif^D r\, \vec{t}(\vec{r})\cdot\left(\nabla_0\cdot\t{P}^T\left(\t{F}(\vec{r})\right)\right)
    =
    -\int_{\Omega_0} \dif^D r\, \left(\nabla_0\otimes\vec{t}(\vec{r})\right) : \t{P}^T\left(\t{F}(\vec{r})\right)
    = 0,
\end{equation}
where $\vec{t}(\vec{r})$ is an arbitrary periodic (vector-valued) test function, $\otimes$ the outer product, $D$ the spatial dimension and $:$ the double dot product, a tensor contraction over two indices, $ \t{A}:\t{B} = A_{ij}B_{ji}$. Surface terms that should appear in Eq.~\eqref{eq:stat-mech-equi-weak-form} vanish due to periodicity. If we interpret the test function $\vec{t}(\vec{r})$ as a displacement, then $\delta \t{F}(\vec{r})=\nabla_0\otimes\vec{t}(\vec{r})$ is a suitable set of compatible test gradients. We can therefore write the equilibrium condition as
\begin{equation}
\begin{split}
    & \int_{\Omega_0} \dif^D r\, \delta\t{F}(\vec{r}) : \t{P}^T\left(\t{F}(\vec{r})\right) \\
    & \qquad =
    \int_{\Omega_0} \dif^D r\, \left(\mathbb{G} \star \delta\tilde{\t{F}}\right)(\vec{r}) : \t{P}^T\left(\t{F}(\vec{r})\right)
    \\
    & \qquad =
    \int_{\Omega_0} \dif^D r\, \delta\tilde{\t{F}}(\vec{r}) : \left(\mathbb{G} \star \t{P}^T\left(\t{F}\right)\right)(\vec{r})
    =
    0
    \label{eq:stat-mech-equi2}
\end{split}
\end{equation}
where now $\delta\tilde{\t{F}}$ is an arbitrary (no longer necessarily compatible) test function and $\mathbb{G} \star \t{A}$ denotes the application of the self-adjoint \emph{operator} $\mathbb{G}$ to a right-hand side object $\t{A}$.
We now discretize the gradients rather than the displacement field, i.e. in the spirit of the Galerkin method, we express the test gradient $\delta\tilde{\t{F}}$ and the deformation gradient $\t{F}$ within the same basis set.

Equation~\eqref{eq:stat-mech-equi2} no longer contains gradients of the constitutive law, given by $\t{P}(\t{F})$. The compatibility operator $\mathbb{G}$ is clearly independent of it since it just ensures fulfillment of the compatibility condition, i.e. $\nabla_0\times (\mathbb{G}:\delta\tilde{\t{F}})=0$ for finite strain formulations. The numerical scheme can therefore be formulated without the explicit requirement of a homogeneous reference material.

The compatibility operator $\mathbb{G}$ is block-diagonal in the Fourier basis, where the blocks are fourth order tensors. This leads to the expression~\cite{Zeman2017}
\begin{equation}
    \label{eq:stat-mech-equi-FFT-homogenization}
    \mathbb{\hat{G}}(\vec{k}):\mathcal{F}\{\t{P}^T(\t{F})\}(\vec{k})=\t{0}
    \quad\text{or}\quad
    \hat{G}_{i \alpha \beta j}(\vec{k}) \mathcal{F}\{P_{j\beta}(\t{F})\}(\vec{k}) = 0,
\end{equation}
where $\mathbb{\hat{G}}$ is  a fourth-order tensor and $\vec{k}$ the wavevector.
We use the hat symbol to denote a Fourier-transformed quantity and $\mathcal{F}\{\cdot\}(\vec{k})$ for the explicit Fourier transform of a quantity as given in \ref{app:discrete-fourier-trafo} by Eq. \eqref{eq:discrete-fourier-trafo}. The numerical solution of Eq.~\eqref{eq:stat-mech-equi-FFT-homogenization} benefits from the fact that for gradient-based optimizers, the steps of the optimization procedure automatically lead to compatible $\t{F}$, see \cite{deGeus2017}.

The operator $\mathbb{G}$ projects arbitrary tensor fields onto \emph{compatible} fields, i.e. fields that can be expressed as the gradient of a lower order tensor. For enforcing the compatibility for second-order tensor fields, the operator is given by \cite{Vondrejc2014}
\begin{equation}
    \label{eq:original-G}
    \hat{G}_{i \alpha \beta j}(\vec{k}) = \delta_{ij} \hat{g}_{\alpha \beta}(\vec{k}), \text{ with }
    \hat{g}_{\alpha \beta}(\vec{k}) = 
    \begin{cases}
        0 & \text{if}\,\vec{k}=\vec{0} , \\
        \frac{k_\alpha k_\beta}{k^2} & \text{else,}
    \end{cases}
\end{equation}
where $k_\alpha$ is the component $\alpha$ of the wavevector $\vec{k}$ and $k=|\vec{k}|$.

\subsection{Interpreting the projection operator}

The projection operator lends itself to a simple interpretation. Let us assume we have an arbitrary vector field $\vec{v}(\vec{r})$. This field is only a gradient $\vec{v}(\vec{r})=\nabla_0 \phi(\vec{r})$ of a scalar field $\phi(\vec{r})$ if its curl vanishes, $\nabla_0\times\vec{v}(\vec{r})=\vec{0}$. If the field is also periodic, we call these fields \emph{compatible}. In the context of homogenization, $\vec{v}$ is a row of the deformation gradient $\t{F}$ and $\phi$ a component of the placement map $\vec{\chi}$. We are interested in the special case where $\vec{v}(\vec{r})$ is periodic but $\phi(\vec{r})$ is not. The periodicity of $\vec{v}$ will be later intrinsically fulfilled through the Fourier transform and we only investigate the gradient property here. For non-compatible fields, we search for the scalar field $\phi(\vec{r})$ that minimizes the residual vector
\begin{equation}
    \vec{R}(\vec{r}) = \nabla_0 \phi(\vec{r}) - \vec{v}(\vec{r})
    \label{eq:residual}
\end{equation}
in a suitable sense (for compatible fields, $\vec{R}\equiv \vec{0}$). For a minimal residual vector ($\vec{R}(\vec{r}) \rightarrow \text{min}$) Equation~\eqref{eq:residual} is an equivalent formulation of Eq.~\eqref{eq:def-grad} and the operator $\mathcal{D}^{-1}$ introduced below is the Green's function of that equation.

The canonical requirement is minimization in the least-squares sense, i.e. minimization of
\begin{equation}
\mathcal{R} = \int_{\Omega_0} d^3 r\, \vec{R}(\vec{r})\cdot\vec{R}(\vec{r}).
\end{equation}
We now need to choose a specific basis set for a series expansion of $\phi(\vec{r})$. In a Fourier basis,
\begin{equation}
    \label{eq:inverse-Fourier-trafo}
    \phi(\vec{r})
    =
    \vec{v}_0\cdot\vec{r}
    + \frac{1}{N}
    \sum_{\vec{k}\not=\vec{0}} \hat{\phi}(\vec{k}) \exp\left(\ii \vec{k}\cdot\vec{r}\right)
\end{equation}
and an equivalent expansion holds for $\vec{v}$,
\begin{equation}
    \label{eq:inverse-Fourier-trafo-v}
    \vec{v}(\vec{r})
    = \frac{1}{N}
    \sum_{\vec{k}} \hat{\vec{v}}(\vec{k}) \exp\left(\ii \vec{k}\cdot\vec{r}\right)
\end{equation}
with $\vec{v}(\vec{0})=\vec{v}_0$ and $N$ the total number of voxels, see \ref{app:discrete-fourier-trafo}. Note that in Eq.~\eqref{eq:inverse-Fourier-trafo} we have set the (arbitrary) mean value of $\phi(\vec{r})$ to zero but added a linear function that cannot be represented in the Fourier basis whose derivative gives the mean value of Eq.~\eqref{eq:inverse-Fourier-trafo-v}. In terms of homogenization, this mean value describes the affine deformation of the whole RVE and plays the role of the boundary condition of the constituting differential equation.

 In the Fourier space the residual becomes
\begin{equation}
    \vec{R}(\vec{k}) = \hat{\vec{\mathcal{D}}}(\vec{k})\hat{\phi}(\vec{k}) - \hat{\vec{v}}(\vec{k})
\end{equation}
for $\vec{k}\not=\vec{0}$ where $\hat{\vec{\mathcal{D}}}(\vec{k})=\ii \vec{k}$ is the Fourier representation of the gradient $\nabla_0$. For $\vec{k}=\vec{0}$ we obtain $\vec{v}_0=\hat{\vec{v}}(0)$.
Parseval's theorem yields $\mathcal{R}=\sum_{\vec{k}} \vec{R}^*(\vec{k})\cdot \vec{R}(\vec{k})$, or
\begin{equation}
    \mathcal{R}
    =
    \sum_{\vec{k}\not=0} \left(
        \hat{\vec{v}}^*\cdot \hat{\vec{v}}
        - \hat{\vec{\mathcal{D}}}\cdot\hat{\vec{v}}^*\hat{\phi}
        - \hat{\vec{\mathcal{D}}}^*\cdot\hat{\vec{v}}\hat{\phi}^*
        + \hat{\vec{\mathcal{D}}}\cdot\hat{\vec{\mathcal{D}}}^* \hat{\phi}^* \hat{\phi}
    \right)
    \label{eq:residualfull}
\end{equation}
where the star is the complex conjugate. Note that all symbols in Eq.~\eqref{eq:residualfull} -- $\hat{\vec{v}}$, $\hat{\vec{\mathcal{D}}}$ and $\hat{\phi}$ -- are function that depend explicit on $\vec{k}$ but that dependence has been omitted here for brevity. Minimization gives the secular equation
\begin{equation}
    \hat{\vec{\mathcal{D}}}(\vec{k})\cdot\hat{\vec{\mathcal{D}}}^*(\vec{k}) \hat{\phi}(\vec{k}) = \hat{\vec{\mathcal{D}}}^*(\vec{k}) \cdot \hat{\vec{v}}(\vec{k}).
\end{equation}
We can solve this for (note that $\vec{k}\not=\vec{0}$)
\begin{equation}
    \hat{\phi}(\vec{k}) = \frac{\hat{\vec{\mathcal{D}}}^*(\vec{k})}{\hat{\vec{\mathcal{D}}}(\vec{k}) \cdot \hat{\vec{\mathcal{D}}}^*(\vec{k})} \cdot \hat{\vec{v}}(\vec{k})
    \equiv \hat{\vec{\mathcal{D}}}^{-1}(\vec{k}) \cdot \hat{\vec{v}}(\vec{k})
\end{equation}
where we interpret the term
\begin{equation}
    \hat{\vec{\mathcal{D}}}^{-1}(\vec{k})
    =
    \frac{\hat{\vec{\mathcal{D}}}^*(\vec{k})}{\hat{\vec{\mathcal{D}}}(\vec{k}) \cdot \hat{\vec{\mathcal{D}}}^*(\vec{k})}
    \label{eq:integration_operator}
\end{equation}
as the inverse of the derivative, i.e. as some form of ``integration''.

In a next step, we compute the gradient of $\hat{\phi}(\vec{k})$. This yields
\begin{equation}
    \hat{\vec{w}}(\vec{k})
    =
    \hat{\vec{\mathcal{D}}}(\vec{k}) \hat{\phi}(\vec{k})
    =
    \hat{\vec{\mathcal{D}}}(\vec{k})
    \left(
    \hat{\vec{\mathcal{D}}}^{-1}(\vec{k})
    \cdot \hat{\vec{v}}(\vec{k})
    \right)
    =
    \hat{\t{g}}(\vec{k}) \cdot \hat{\vec{v}}(\vec{k})
    \label{eq:fourierprojection}
\end{equation}
with
\begin{equation}
    \hat{\t{g}}(\vec{k})
    =
    \hat{\vec{\mathcal{D}}}(\vec{k})\otimes\hat{\vec{\mathcal{D}}}^{-1}(\vec{k}).    
\end{equation}
The operator $\hat{\t{g}}(\vec{k})$ hence projects an arbitrary field on its compatible part in the least squares sense with respect to the integral inner product ($L^2$-norm). The full projection operator for the deformation gradient is given by the first part of Eq.~\eqref{eq:original-G}. Using the Fourier derivative $\hat{\vec{\mathcal{D}}}(\vec{k})=\ii \vec{k}$ yields the specific form of the projection operator given in Eq.~\eqref{eq:original-G}. Since $\mathbb{G}$ contains the Green's function of Eq.~\eqref{eq:def-grad}, but \emph{not} of Eq.~\eqref{eq:stat-mech-equi}, it is independent of the constitutive law and does not require a reference material. The formulation of the projection operator for small-strain elasticity is described in \ref{sec:small-strain}.

\subsection{Discrete projection}

Rather than using Eq.~\eqref{eq:inverse-Fourier-trafo},
we can expand $\phi(\vec{r})$ in other bases of choice. For example, we will below employ linear finite elements.
Other discretizations of the gradient operator $\nabla_0$ with less suitable properties can be obtained through finite-differences methods. We here assume that the simulation cell is structured in a regular (equally spaced) grid with node positions $\{\vec{r}_{IJ}\}$, where $I$ and $J$ are node indices. We will call the individual grid cell a \emph{voxel} and develop the theory in two dimensions, but generalization to three dimensions is straightforward. The placement map $\vec{\chi}$ is only known at the nodes which are the corners of the voxels.

The discrete derivative (in some ``direction'' $\alpha$) in voxel $I,J$ can generally be written as the convolution
\begin{equation}
    \mathcal{D}_\alpha \chi(\vec{r}_{IJ}) = \frac{1}{\Delta^{(\alpha)}} \sum_{ij} s^{(\alpha)}_{ij} \chi(\vec{r}_{I+i,J+j})
    \label{eq:discrderiv}
\end{equation}
with $\vec{r}_{I+i,J+j}=\vec{r}_{IJ}+\vec{r}_{ij}$ and $\Delta^{(\alpha)}$ the voxel size in direction $\alpha$ where the brackets and the upper index indicate that there is no Einstein sum convention applied for the index $(\alpha)$ and periodicity of $\vec{r}$ in a natural sense is considered.
The collection of coefficients $s_{ij}$ is called the \emph{stencil} of the operation. We will introduce different stencils below and note that derivatives in different directions $\alpha$ require different stencils. In the following we will assume that the deformation gradient is described by a set of $d$ derivatives, but that $d$ is not necessarily equal to the dimension $D$ of the system. This allows the subdivision of voxels into multiple evaluation points, i.e., $d$ may be an integer multiple of $D$. Using $n_q$ evaluation points per voxel leads to $d=n_q D$ derivatives per voxel within the framework that we now develop. The subscript $()_q$ was chosen because these evaluation points correspond to Gaussian quadrature points in a classic finite element discretization.

By expanding the discrete placement map into a Fourier series of the form given by Eq.~\eqref{eq:inverse-discrete-fourier-trafo}, we can write the derivative operation in Fourier space as
\begin{equation}
    \hat{\mathcal{D}}_\alpha(\vec{k}) = \frac{1}{\Delta^{(\alpha)}}\sum_{ij} s^{(\alpha)}_{ij} \exp\left(\ii\vec{k}\cdot\vec{r}_{ij} \right).
    \label{eq:genfourderiv}
\end{equation}
Again the case $\vec{k}=\vec{0}$ is special since $\sum_{ij} s^{(\alpha)}_{ij}=0$.
For $\vec{k}\not=\vec{0}$, we can use the same argument as above: A general vector field $\hat{\vec{v}}(\vec{k})$ can be projected (in the least squares sense) onto its compatible part $\hat{\vec{w}}(\vec{k})$ using
\begin{equation}
    \hat{w}_\alpha(\vec{k}) = \hat{g}_{\alpha\beta}(\vec{k}) \hat{v}_\beta(\vec{k})
\end{equation}
with
\begin{equation}
\label{eq:discrete_projection_op}
    \hat{g}_{\alpha\beta}(\vec{k}) = \frac{\hat{\mathcal{D}}_\alpha(\vec{k}) \hat{\mathcal{D}}^*_\beta(\vec{k})}{\hat{\mathcal{D}}_\gamma(\vec{k})\hat{\mathcal{D}}^*_\gamma(\vec{k})} 
    = \hat{\mathcal{D}}_\alpha(\vec{k}) \hat{\mathcal{D}}^{-1}_\beta(\vec{k}),
\end{equation}
cf. Eq. \eqref{eq:integration_operator}. We want to mention that by the Helmholtz decomposition $\delta_{\alpha\beta}-g_{\alpha\beta}(\vec{k})$ is a projection to divergence free fields, with applications to error estimation. The projection operator for the deformation gradient is then given by $\hat{G}_{i\alpha\beta j}(\vec{k}) = \delta_{ij} \hat{g}_{\alpha\beta}(\vec{k})$. Note that $\hat{\t{g}}\in \mathbb{C}^{d\times d}$ and $\hat{\mathbb{G}}\in \mathbb{C}^{D\times d\times d\times D}$. The projection operator is idempotent (i.e. a \emph{projection}), since,
\begin{equation}
\begin{split}
    \hat{G}_{l \gamma \alpha i}(\vec{k}) & \hat{G}_{i \alpha \beta j}(\vec{k}) \hat{f}_{j \beta}(\vec{k}) \\
    &= \delta_{li} \hat{g}_{\gamma \alpha}(\vec{k}) \left( \delta_{ij} \hat{g}_{\alpha \beta}(\vec{k}) \hat{f}_{j \beta}(\vec{k}) \right) \\
    &= \delta_{li} \delta_{ij} \hat{\mathcal{D}}_\gamma(\vec{k}) \hat{\mathcal{D}}^{-1}_\alpha(\vec{k}) \hat{\mathcal{D}}_\alpha(\vec{k}) \hat{\mathcal{D}}^{-1}_\beta(\vec{k}) \hat{f}_{j \beta}(\vec{k}) \\
    &= \delta_{lj} \hat{g}_{\gamma \beta}(\vec{k}) \hat{f}_{j \beta}(\vec{k}) \\
    &= \hat{G}_{l \gamma \beta j}(\vec{k}) \hat{f}_{j \beta}(\vec{k})
\end{split}
\end{equation}
This compatibility operator $\hat{\mathbb{G}}$ is the projection generalized for arbitrary derivative operators $\hat{\vec{\mathcal{D}}}$.

The case $\vec{k}=\vec{0}$ contains the boundary condition and is treated as for the Fourier derivative. For multiple elements, each element needs to hold the same gradient for $\vec{k}=\vec{0}$. We want to emphasize that in this discrete basis, the discrete Fourier transformation has the role of accelerating the convolution rather than providing the basis set for the underlying discretization.

\subsection{Finite differences and Lanczos-$\sigma$ correction}

Canonical discrete derivative operators are given by finite difference schemes. We here consider the first-order central differences,
\begin{equation}
    \frac{\partial\chi_i}{\partial r_\alpha} = \frac{\chi_i\left(r_j + \delta_{j\alpha} \Delta^{(\alpha)}\right) - \chi_i\left(r_j - \delta_{j\alpha} \Delta^{(\alpha)}\right)}{2 \Delta^{(\alpha)}} + \mathcal{O}\left((\Delta^{(\alpha)})^2\right),
    \label{eq:central-difference-derivative-real}
\end{equation}
which can be expressed in Fourier space as follows (see e.g.~\cite{Vidyasagar2017}): \begin{equation}
    \hat{\mathcal{D}}_\alpha^{\text{cd}}(\vec{k}) = \frac{\ii\sin(k_\alpha \Delta^{(\alpha)})}{\Delta^{(\alpha)}}.
    \label{eq:central-difference-derivative-fourier}
\end{equation}
Remember that $\Delta^{(\alpha)}$ is as before the grid spacing in direction $\alpha$ and there is no summation over upper indices in parenthesis. 
We note that this central differences scheme is related to the Lanczos-$\sigma$ correction for Gibbs ringing in direction $\alpha$~\cite{Hamming1962}, $\sigma^{(\alpha)}(\vec{k})=\sin(k_\alpha\Delta^{(\alpha)})/k_\alpha\Delta^{(\alpha)}$. Expressed using the $\sigma$-factor, the derivative operator is given by $ \mathcal{D}^{\text{cd}}_\alpha(\vec{k}) = \ii k_\alpha \sigma^{(\alpha)}(\vec{k})$.

\begin{figure*}
    \centering
    \includegraphics[width=\textwidth]{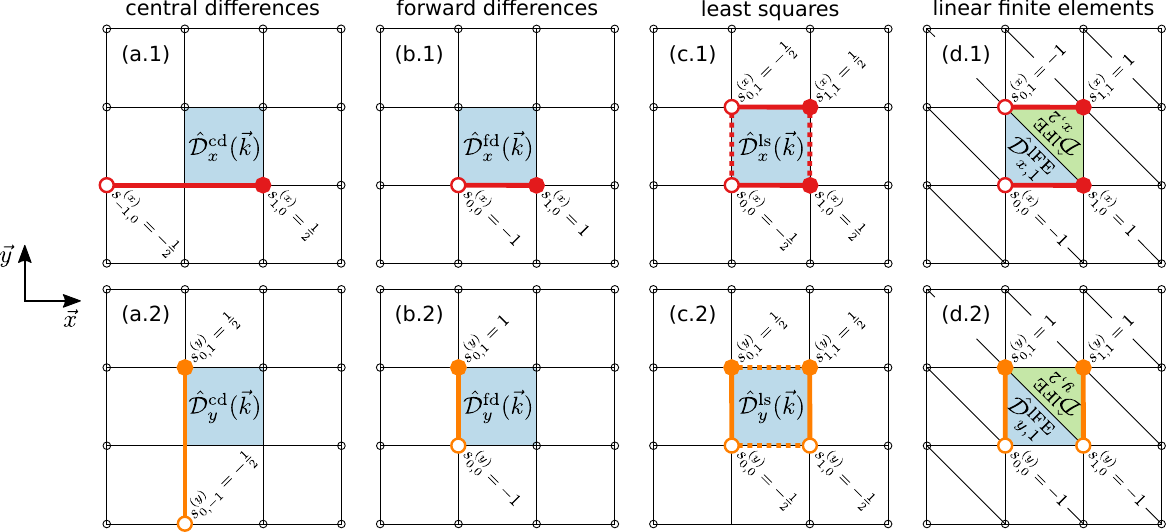}
    \caption{Graphical representation of the different stencils $s^{(x)}_{ij}$ and $s^{(y)}_{ij}$ for the derivatives in $x$ and $y$-direction of the discrete derivative operators $\hat{\mathcal{D}}_\alpha(\vec{k})$. Column \textbf{(a)} shows central differences, \textbf{(b)} forward differences, \textbf{(c)} least squares and \textbf{(d)} linear finite element stencils. Row \textbf{(1)} shows the derivatives in $x$-direction and \textbf{(2)} in $y$-direction. Computed derivatives are assigned to the voxel marked in blue. For linear finite elements in \textbf{(d)}, the voxel is subdivided into two triangles. Full dots indicate positive stencil values and open dots negative ones. Thick lines indicate the direction of the derivative and dotted lines in column \textbf{(c)} indicate connected stencils in the non derivative direction.}
    \label{fig:stencils}
\end{figure*}

Another common first-order finite difference scheme is forward differences, 
\begin{equation}
    \frac{\partial\chi_i}{\partial r_\alpha} = \frac{\chi_i(r_j + \delta_{j\alpha} \Delta^{(\alpha)}) - \chi_i(r_j)}{\Delta^{(\alpha)}} + \mathcal{O}(\Delta^{(\alpha)})
\end{equation}
with the Fourier-space representation
\begin{equation}
    \hat{\mathcal{D}}_\alpha^{\text{fd}}(\vec{k}) = \frac{\exp\left(\ii k_\alpha \Delta^{(\alpha)}\right) - 1}{\Delta^{(\alpha)}}.
    \label{eq:forward-difference-derivative-fourier}
\end{equation}
The full stencil coefficients for these two schemes as applied to two-dimensional problems are shown in Fig.~\ref{fig:stencils}a and \ref{fig:stencils}b. Inserting these coefficients into the generic expression Eq.~\eqref{eq:genfourderiv} yields the specific Fourier representations given in Eqs.~\eqref{eq:central-difference-derivative-fourier} and \eqref{eq:forward-difference-derivative-fourier}. Figure~\ref{fig:stencils} also has a graphical representation of these discrete derivatives.

\subsection{Least squares}
\label{sec:least-squares}

\begin{figure}
    \centering
    \includegraphics[width=0.5\textwidth]{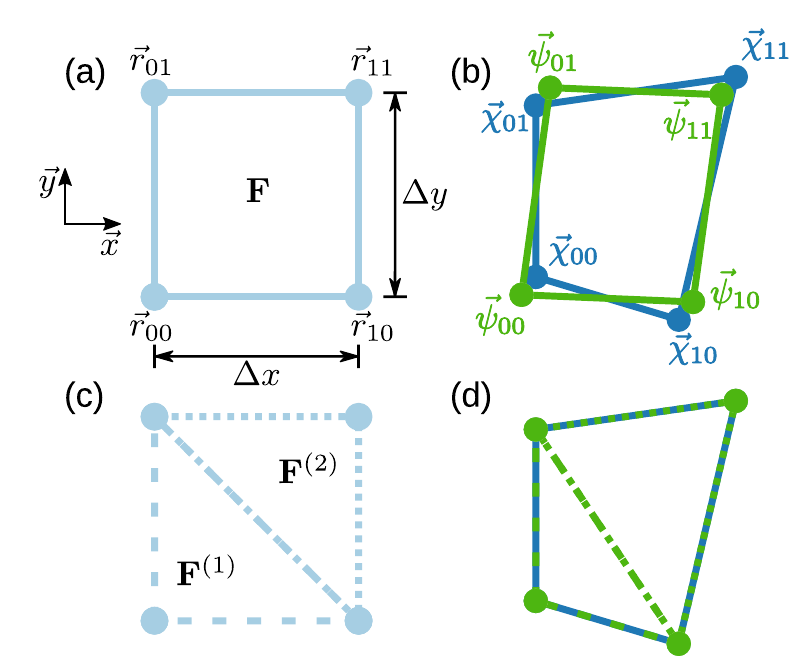}
    \caption{Non-affine deformation of a single 2D voxel. \textbf{(a)} The undeformed voxel (light blue) is described by its four corners $\vec{r}_{ij}$ in a grid with grid spacing $\Delta x$ and $\Delta y$. \textbf{(b)} A non-affine deformation of this voxel is shown in dark blue. The green parallelogram shows the least-squares approximation of the dark blue trapezoid using a uniform gradient $\t{F}$ througout the voxel. \textbf{(c)} The single voxel of \textbf{(a)} is subdivided into two triangles, the deformation of each is described by a uniform deformation gradients $\t{F}^{(1)}$ for the dashed triangle and $\t{F}^{(2)}$ for the dotted triangle. \textbf{(d)} Individual affine deformation of these two triangles (green) can perfectly match the prescribed corner displacements (dark blue).}
    \label{fig:single-element}
\end{figure}

We now turn to a purely geometric interpretation of the deformation gradient to derive alternative discrete stencils. The deformation maps an infinitesimal fiber vector $\delta\vec{r}$ into $\delta\vec{r}'=\t{F(}\vec{r})\cdot\delta\vec{r}$. Assuming constant $\t{F}$ over a given voxel (light blue rectangle in Fig.~\ref{fig:single-element}a), the voxel is affinely deformed by the deformation gradient $\t{F}$ (green parallelogram in Fig.~\ref{fig:single-element}b) and cannot represent arbitrary displacements of the corners (dark blue trapezoid in Fig.~\ref{fig:single-element}b) as long as the deformation gradient $\t{F}$ is uniform on that voxel. In order to represent the corner displacements exactly, we can for instance subdivide the voxel, e.g. in 2D into two triangles (see Fig.~\ref{fig:single-element}c and Fig.~\ref{fig:single-element}d), introducing multiple elements per voxel, with their own uniform deformation gradient per element. We will discuss this decomposition in the next section and will for now continue to work with a uniform deformation gradient per voxel and require matching of the corner displacements in a least-squares sense.

We now regard the displaced corner positions only within a single voxel (see Fig.~\ref{fig:single-element}a). The undeformed coordinates of the voxel are $\vec{r}_{kl}$ with $k,l\in\{0,1\}$. The true deformed coordinates are $\vec{\chi}_{kl}$ and the affinely deformed coordinates produced by the deformation gradient $\t{F}$ are $\vec{\psi}_{kl}$. We require that the affinely deformed rectangle matches the true deformation in a least square sense. This means we are looking for the deformation gradient $\t{F}$ that minimizes the residual
\begin{equation}
    R = \sum_{k,l} \left| \vec{\psi}_{kl} - \vec{\chi}_{kl} \right|^2
    \label{eq:lstsqres}
\end{equation}
Since the affinely deformed voxel is spanned by the two vectors $\t{F}\cdot\Delta \vec{r}_{10}$ and $\t{F}\cdot\Delta\vec{r}_{01}$ with $\Delta \vec{r}_{kl}=\left(\vec{r}_{kl} - \vec{r}_{00}\right)$, we can write the affinely deformed coordinates as
\begin{align}
    \vec{\psi}_{10} &= \vec{\psi}_{00} + \t{F}\cdot \Delta\vec{r}_{10} \\ 
    \vec{\psi}_{01} &= \vec{\psi}_{00} + \t{F}\cdot \Delta\vec{r}_{01} \\ 
    \vec{\psi}_{11} &= \vec{\psi}_{00} + \t{F}\cdot \Delta\vec{r}_{11}.
\end{align}
We now insert these into Eq.~\eqref{eq:lstsqres} and minimize with respect to $\t{F}$ and $\vec{\psi}_{00}$.
This yields the secular equations
\begin{align}
    \frac{\partial R}{\partial \vec{\psi}_{00}}
    &=
    2\sum_{k,l} \left(\vec{\psi}_{kl} - \vec{\chi}_{kl}\right)
    =
    \vec{0}
    \\
    \frac{\partial R}{\partial \t{F}}
    &=
    2\sum_{k,l} \left(\vec{\psi}_{kl} - \vec{\chi}_{kl}\right) \otimes \Delta\vec{r}_{kl}
    =
    \t{0}
\end{align}
or
\begin{equation}
    \vec{\psi}_{00}
    =
    \frac{1}{4}
    \sum_{k,l} \left(\vec{\chi}_{kl} - \t{F}\cdot \Delta\vec{r}_{kl}\right)
\end{equation}
and
\begin{equation}
    \t{F}\cdot\sum_{k,l}\left(\Delta\vec{r}_{kl}-\frac{1}{4}\sum_{m,n}\Delta\vec{r}_{mn}\right)\otimes\Delta\vec{r}_{kl}
    =
    \sum_{k,l} \left(\vec{\chi}_{kl} - \frac{1}{4}\sum_{m,n}\vec{\chi}_{mn}\right)\otimes\Delta\vec{r}_{kl}.
\end{equation}
For rectangular lattices with lattice spacing $\Delta x$ and $\Delta y$, this can be solved to give
\begin{equation}
    \t{F}
    =
    \frac{1}{2}
    \begin{pmatrix}
        \frac{\vec{\chi}_{10} - \vec{\chi}_{00} + \vec{\chi}_{11}-\vec{\chi}_{01}}{\Delta x} & \frac{\vec{\chi}_{01} - \vec{\chi}_{00} + \vec{\chi}_{11}-\vec{\chi}_{10}}{\Delta y}
    \end{pmatrix}.
\end{equation}
The corresponding stencil coefficients are shown in Fig.~\ref{fig:stencils}c.

\subsection{Linear finite elements}
\label{sec:linear_finite_elements}

The previous section argued that a uniform deformation gradient per voxel is insufficient to represent the voxel's deformation as measured by the displacement of each corner. This becomes evident from simply counting the degrees of freedom: In two dimensions, the voxel's deformation (and rotation) is given by three vectors ($6$ degrees of freedom) while the deformation gradient has $4$ independent components. In this section, the problem is solved by splitting the voxel into two triangles that are each described by a uniform deformation gradient (see Fig.~\ref{fig:single-element}c). The voxel still has 6 degrees of freedom, but now we have 2 deformation gradients with a total of 8 independent components and the constraint that the diagonal boundary between the two triangles remain of the same length and direction (two blocked degrees of freedom). From a simple geometric argument, this triangular decomposition can hence describe arbitrary corner displacements (see Fig.~\ref{fig:single-element}d).

This geometric point of view is fully equivalent to a formulation using linear finite elements. Within our rectangular lattice, we use the usual linear shape functions
\begin{align}
    N^{(1)}_{00}(x,y) &= 1-x/\Delta x-y/\Delta y \\
    N^{(1)}_{10}(x,y) &= x/\Delta x \\
    N^{(1)}_{01}(x,y) &= y/\Delta y
\end{align}
where the origin of the coordinate system is at the bottom left of the voxel and an equivalent set of shape functions exists for element $(2)$ in Fig.~\ref{fig:single-element}c. The shape function gradient of $\chi(x,y)=\chi_{00} N^{(e)}_{00}(x,y)+\chi_{10} N^{(e)}_{10}(x,y)+\chi_{01} N^{(e)}_{01}(x,y)$ is constant on element $(e)$ and given by the stencil coefficients shown in Fig.~\ref{fig:stencils}d for the two deformation gradients. The stencil for element $(1)$ is identical to the forward differences scheme. The stencil for element $(2)$ is a forward differences scheme evaluated on a different set of nodes turning it into a backward differences scheme. Generalizations of this scheme to non-orthogonal voxels and three dimensions are straightforward. This formulation is identical to traditional linear finite elements, but unlike classical finite-elements, the projection formulation yields a condition number that is independent of system size, leading to scale-independent convergence properties~\cite{Ladecky2021}. Note that the least square approach described in the previous section simply yields the average of the deformation gradients on the two triangles.

It is important to emphasize that this formulation requires storing two deformation gradients per voxel. In the discrete projection developed above, this means we have $n_q = 2$ evaluation points and $d=4$ derivatives for the two-dimensional formulation. The deformation gradient $\t{F}$ and Piola-Kirchhoff stress $\t{P}$ are then both elements of $\mathbb{R}^{D\times D \times n_q}$ within the projection scheme. For the evaluation of the constitutive law, both $\t{F}$ and $\t{P}$ need to be decomposed in their element-wise contributions $\t{F}^{(e)}$ and $\t{P}^{(e)}$ that are elements of $ \mathbb{R}^{D\times D}$. We can formally write
\begin{equation}
    \t{F} = \begin{pmatrix}
        \t{F}^{(1)} \\
        \t{F}^{(2)}
    \end{pmatrix}
\end{equation}
for this decomposition. We note that finite element discretizations with multiple quadrature points can be described in this discrete projection using similar decompositions.

\begin{figure}
    \centering
    \includegraphics[width=0.7\textwidth]{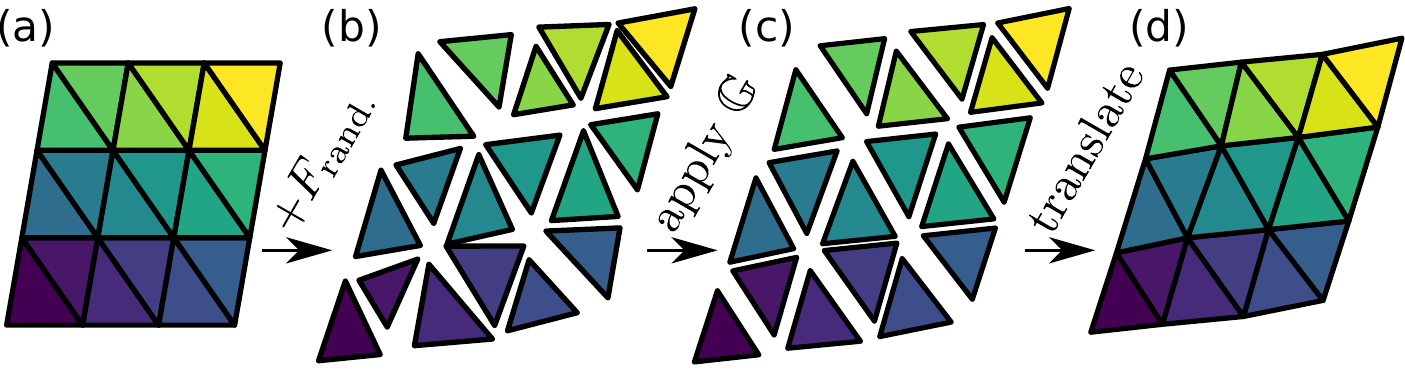}
    \caption{The suggested projection operator for the finite-element formulation is acting on a random, non-compatible deformation gradient field and projects it onto the respective compatible field. The figure shows \textbf{(a)} an undeformed grid \textbf{(b)} in which triangular elements are independently randomly deformed and presented in an exploded view. The random deformation is corrected by \textbf{(c)} applying the projection operator $\mathbb{G}$ and \textbf{(d)} translating from the exploded view into a compact cell where all triangles fit together.}
    \label{fig:effect_of_G}
\end{figure}

The interpretation put forward at the beginning of this section is that the deformation gradient describes the geometry of the respective triangle (see Fig.~\ref{fig:single-element}d). This allows an intuitive interpretation of the projection operator $\mathbb{G}$. We start with a decomposition of some domain into triangles (Fig.~\ref{fig:effect_of_G}a). The rotation and shape of each of these triangles is described by a deformation gradient $\t{F}^{(e)}\in \mathbb{R}^{D\times D}$. We now randomly disturb these triangles by adding a random number to the components of their deformation gradients. The resulting structure (see exploded view in Fig.~\ref{fig:effect_of_G}b) is clearly no longer compatible. Application of $\mathbb{G}$ (Fig.~\ref{fig:effect_of_G}c) adjusts the shape of each triangle such that they are again compatible and can be assembled into a continuous deformed structure (Fig.~\ref{fig:effect_of_G}d).

\subsection{Even vs. odd number of grid points}

The original formulation of the compatibility projection that employs the Fourier-derivative only works exactly for odd-sized grids, which is a result of the structure of the projection operator. The Fourier derivative is ambiguous at the Nyquist frequency (or the edge of the first Brillouin zone). This ambiguity originates in the freedom of choosing one of two possible equivalent even numbered Fourier grids $\left\{k_i\right\}$ from the class $k_i \Delta \in \left[-\pi,\pi\right)$ or $k_i 
\Delta \in \left(-\pi,\pi\right]$, where $\Delta$ is the grid spacing. Even sized grids sample the frequency exactly at the Brillouin zone boundary (see Eq. \eqref{eq:Fourier_wave_vectors} for the choice $k_i \Delta \in \left[-\pi,\pi\right)$). As a result the two possible grids differ only in one single grid point, the Nyquist frequency.

\begin{figure}
    \centering
    \hspace{-0.4cm}
    \includegraphics[width=0.7\columnwidth]{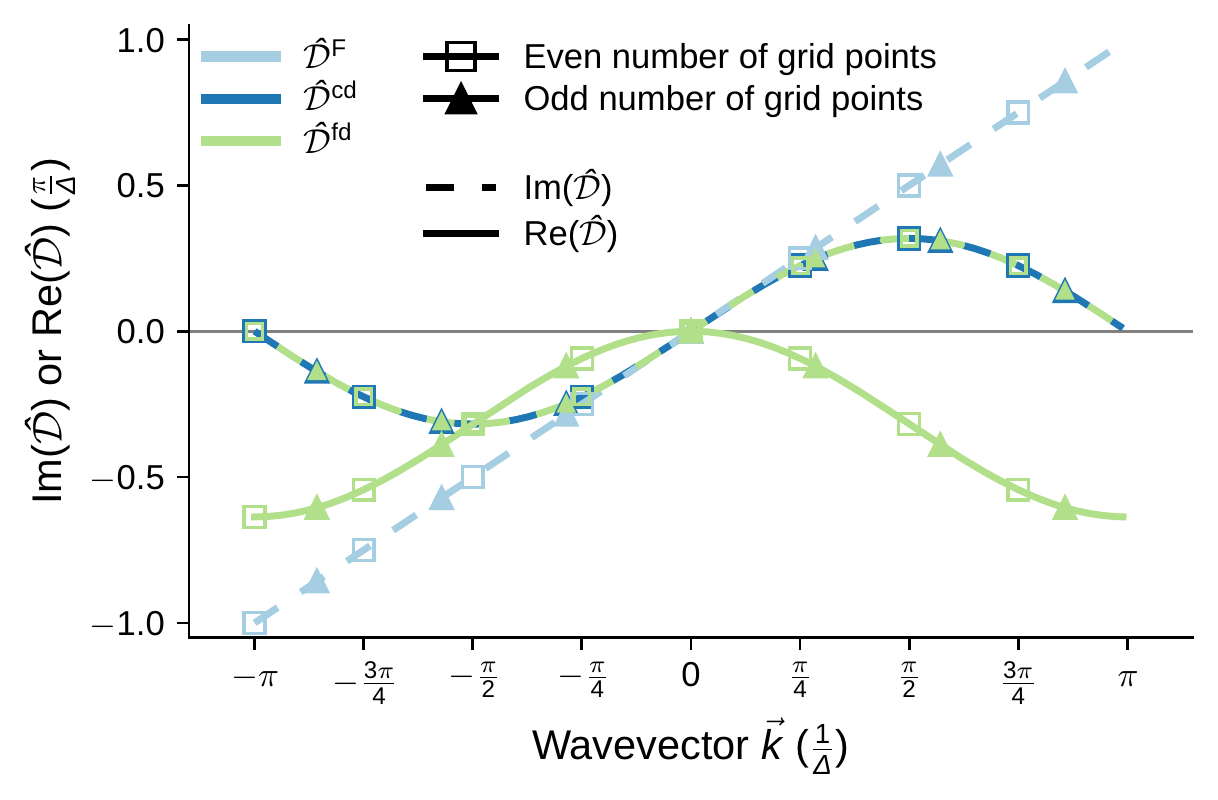}
    \vspace{-0.3cm}
    \caption{Real (continuous line) and imaginary (dashed line) parts of the three different derivative operators, $\hat{\mathcal{D}}^\text{F}$ (Fourier derivative, light blue), $\hat{\mathcal{D}}^\text{cd}$ (central differences, dark blue) and $\hat{\mathcal{D}}^\text{fd}$ (forward differences, green). The real part of the central differences and Fourier-type projection operator are zero everywhere and therefore not shown. Data points are shown for an even grid (open squares) with eight points and for an odd grid (triangles) with seven points, see Eq. \eqref{eq:Fourier_wave_vectors}. Only the even grid has a data point at the Nyquist frequency, here chosen to be at $-\pi/\Delta$. Only the forward differences derivative operator satisfies $\hat{\mathcal{D}}(k_\text{Ny}) \hat{\mathcal{D}}^*(k_\text{Ny}) \neq 0$ and $\hat{\mathcal{D}}(k_\text{Ny}) = \hat{\mathcal{D}}(-k_\text{Ny})$ which are the two requirements for an unambiguous projection operator on arbitrary sized grids.}
    \label{fig:even-odd_projection_operator}
\end{figure}

The Nyquist frequency is given by $k_{\text{Ny}}=\pm \pi/\Delta$, where an even sized grid contains either the positive or negative Nyquist frequency. Typically this small difference does not matter for the periodic function $\hat{f}(k)$ because the Fourier coefficients are also equivalent in the single point $\pm k_\text{Ny}$ where the two grids differ, $\hat{f}(-k_\text{Ny})$ from the one grid containing $-k_\text{Ny}$ has the same value as $\hat{f}(+k_\text{Ny})$ from the other grid containing $+k_\text{Ny}$.
Thus, the Fourier series does not depend whether the positive or negative Nyquist frequency is included in the grid and is therefore unambiguous for any (sampled) function $f(x)$.
However, when analyzing the Fourier series of $\mathcal{D}^\text{F} f(x)$, where the upper index ``$\text{F}$'' indicates the explicit use of the Fourier derivative $\hat{\mathcal{D}}^\text{F}(k)=\ii k$, evaluated on grid points $x$,
\begin{equation}
    \label{eq:Df(x)_Fourier_series}
    \mathcal{D}^\text{F} f(x) = \frac{1}{N} \sum_k \ii k \hat{f}(k) \exp \left(\ii k x\right) = \frac{1}{N} \sum_k S(k,x),
\end{equation}
we find for the summand $S(k_\text{Ny}, x)$ at the Nyquist frequency
\begin{equation}
    S(k_\text{Ny}, x) = \ii k_\text{Ny} \hat{f}(k_\text{Ny}) \exp \left(\ii k_\text{Ny} x\right) = \ii k_\text{Ny} \hat{f}(k_\text{Ny}) \cos (k_\text{Ny} x)
\end{equation}
where we used $k_\text{Ny} x \propto \pi n, n\in\mathbb{Z}$. Since $\hat{\mathcal{D}}^\text{F}(k_\text{Ny}) = \ii k_\text{Ny} \not = - \ii k_\text{Ny} = \hat{\mathcal{D}}^\text{F}(-k_\text{Ny})$ the summand $S(k_\text{Ny}, x) \neq S(-k_\text{Ny}, x)$ which gives an ambiguity in the Fourier series of $\mathcal{D}^\text{F} f(x)$ for even sized grids at the Nyquist frequency. The outcome of the Fourier series in Eq. \eqref{eq:Df(x)_Fourier_series} depends on which one of the two possible even sized grids is taken. This ambiguity vanishes for odd sized grids, since the function is never evaluated at the Nyquist frequency and there is only a single possible choice for the Fourier grid. Figure~\ref{fig:even-odd_projection_operator} plots $\hat{\mathcal{D}}^\text{F}(k)$ at the discrete evaluation points for even and odd-sized grids as an illustration of this problem.

This ambiguity is resolved for all the discrete stencils. For example, the forward differences stencils fulfills $\hat{\mathcal{D}}^\text{fd}(k_\text{Ny})=\hat{\mathcal{D}}^\text{fd}(-k_\text{Ny})$ (see Fig.~\ref{fig:even-odd_projection_operator}).
Since we are interested in the compatibility projection operator we additionally require a second condition $\hat{\mathcal{D}}(k_\text{Ny}) \hat{\mathcal{D}}^*(k_\text{Ny}) \neq 0$ for a proper derivative operator on even grids to prevent a division by zero in the computation of the compatibility operator by Eq.~\eqref{eq:discrete_projection_op}, or in other words we need a properly defined inverse derivative operator $\hat{\mathcal{D}}^{-1}(k)$ (see Eq.~\eqref{eq:integration_operator}). 
An exception from the presented discrete stencils is the central differences stencil, where $\hat{\mathcal{D}}(k_\text{Ny})\hat{\mathcal{D}}^*(k_\text{Ny})=0$ as a result of a decoupling of the domain into two subgrids. An alternative point of view is that the discrete Fourier-transform is simply employed to carry out a discrete convolution and there is no ambiguity with respect to the Nyquist frequency.

\section{Examples and validation}
\label{chap:examples}

In the following, we describe four two-dimensional examples to demonstrate the methods developed above with a focus on ringing phenomena.
First, we investigate a single soft voxel in a uniform hard matrix under biaxial strain. This minimal example already shows strong ringing artifacts in the $xy$-component of the stress tensor.
Second, we analyse a cell with two pillars separated by one layer of voxels with the Young modulus set to zero. This example shows the ability to handle infinite material contrast with vacuum~\cite{Schneider2020} and simulate a free surface despite the intrinsic periodic boundary conditions~\footnote{Note that simulating infinite phase contrast with rigid pixels is not possible without additional treatment, as it would introduce multiplications with infinity in Eq. \eqref{eq:stat-mech-equi-FFT-homogenization}.}. Additionally, one of the pillars contains an inhomogeneity which gives rise to ringing artifacts in the original spectral formulation. Third, we test the numeric correctness of the results by investigating an Eshelby inhomogeneity. We compare the results obtained by the FFT-based method against the analytical Eshelby solution, corrected for periodic boundary conditions. Finally, we demonstrate the feasibility of the method for complex constitutive laws on a damage mechanics problem.

\subsection{Single voxel inhomogeneity}
\label{sec:single_pixel_inhomo}

\begin{figure*}
    \centering
    \includegraphics[width=\textwidth]{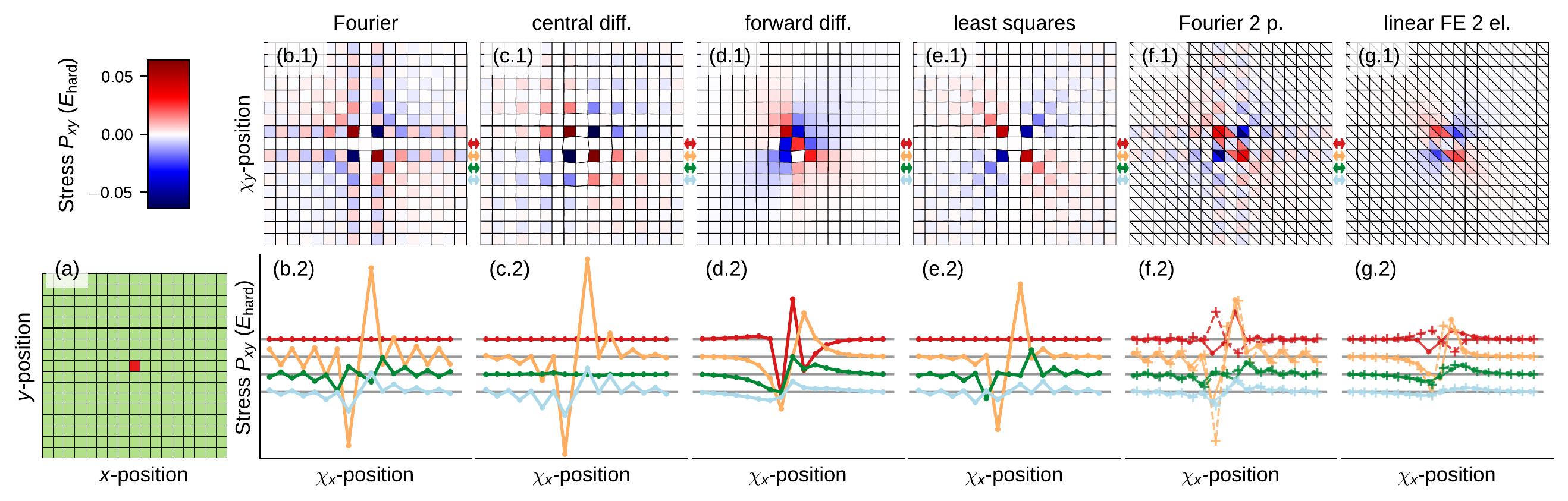}
    \caption{Soft single voxel inhomogeneity in hard matrix ($E_\text{hard}=10 \, E_\text{soft}$) at $10\%$ biaxial strain computed on a $17\times17$ grid. \textbf{(a)} shows the phase setup of the central inhomogeneity (red) embedded in the matrix (green). The first row \textbf{(x.1)} shows the color coded shear component of the first Piola-Kirchhoff stress $P_{xy}$ on the deformed grid. The second row gives the shear stress along rows of the grid as indicated by the colored double arrows in between the subfigures of the first row. In red the row through the inhomogeneity voxel and then going down row by row in orange, green and blue. The point markers indicate the voxel data and the lines are only added to guide the eye. In the last two columns the continuous line with point markers represents the shear stress of the lower triangle and the dashed line with plus markers represents the upper triangle (see Fig.~\ref{fig:single-element}). The data points are presented at the geometric center of the triangles and voxels. The gray line is the level of zero stress for each row and the y-scaling is the same for all columns \textbf{(x.2)} to make a direct comparison possible.
    The columns represent the data generated with different projection operators as indicated by the column titles: \textbf{(b.i)} Fourier, \textbf{(c.i)} central differences, \textbf{(d.i)} forward differences, \textbf{(e.i)} least squares, \textbf{(f.i)} Fourier type on two evaluation points per voxel and \textbf{(g.i)} the linear finite element type projection on two elements per voxel.}
    \label{fig:single_voxel_inhomo}
\end{figure*}

A classical continuum mechanics problem is the inhomogeneity, an inclusion of a material in a matrix with different material properties. At the boundary of the inhomogeneity, there is a discrete change in the material properties which usually leads to Gibbs ringing artifacts in spectral methods. As minimal example of such an inhomogeneity, we present a single voxel inhomogeneity (in red) placed in the center of a $17\times17$ voxel matrix (in green) as shown in Fig.~\ref{fig:single_voxel_inhomo}a. The matrix with Young modulus $E_\text{hard}$ is ten times harder than the inhomogeneity $E_\text{soft}=E_\text{hard}/10$, and both have the same Poisson ratio of $\nu=0.33$. The material response is described by an isotropic finite strain linear elastic law as described in Refs.~\cite{Belytschko2014,deGeus2017}. 

We apply a biaxial tensile strain of 10\% ($F_{xx}=F_{yy}=0.1$ and $F_{xy}=F_{yx}=0$) and investigate the shear component $P_{xy}$ of the first Piola-Kirchhoff stress for different implementations of the projection operator. Figure~\ref{fig:single_voxel_inhomo} gives an overview of the results where the first row, panels (b.1) to (g.1), show the color coded stress and the second row, panels (b.2) to (g.2), give a more detailed look on the stress along selected rows of the matrix as indicated by the colored arrows in between the subfigures in the first row of Fig. \ref{fig:single_voxel_inhomo}. The lines in the panels (b.2) to (g.2) are ordered from top to bottom starting from the center row, i.e. the row with the inhomogeneity in red. 
Each column represents the result found with a different projection operator: Column (b) for the Fourier-type projection operator as given by Eq.~\eqref{eq:original-G}, (c) for central differences given by Eq.~\eqref{eq:central-difference-derivative-fourier}, (d) for forward differences given by Eq.~\eqref{eq:forward-difference-derivative-fourier}, (e) for the least square scheme described in Sec.~\ref{sec:least-squares}, (f) for the Fourier-type projection operator on two evaluation points per voxel (see \ref{app:fourier-G-2-points-per-voxel}) and (g) for linear finite elements on two elements per voxel as derived in Sec.~\ref{sec:linear_finite_elements}.

Panel (b) show the stress field for the original method (e.g. Ref.~\cite{deGeus2017}), with a projection operator based on the Fourier derivative. As expected we observe strong ringing artifacts leading to a checkerboard pattern of the stress field. The stress field and its oscillations are strongest at the inhomogeneity and decay with increasing distance to the discontinuity. However, the ringing does not disappear even at the edges of the cell which was also tested for finer grids, different material contrasts and a slightly inhomogeneous matrix (results not shown). The symmetry of the setup leads to a line of zero stress in the row and column that contains the inhomogeneity (see Fig.~\ref{fig:single_voxel_inhomo}b.1 and red line in Fig.~\ref{fig:single_voxel_inhomo}b.2).

Results obtained with the central differences projection operator are shown in Fig.~\ref{fig:single_voxel_inhomo}c.1 and c.2. The Gibbs ringing artifacts should be strongly suppressed for this method, however we observe a checkerboard pattern of different style compared to (b). This checkerboard pattern originates in the well-known~\cite{book-Ferziger2002,Rauwoens2007} (odd-even) decoupling into two subgrids over short distances of the central differences stencil shown in Fig.\ref{fig:stencils}a. The decoupling of the two grids is not complete because of an odd-sized grid $(17\times17)$ and the oscillations decay with increasing distance from the inhomogeneity.

The forward differences stencil, results shown in panels (d), leads to an oscillation-free but slightly asymmetric stress field which originates from the asymmetry of the forward differences scheme (see Fig.\ref{fig:stencils}b).
This asymmetry is corrected by the least square stencil shown in panels (e). However, the stress has also a checkerboard pattern with checkerboard characteristics of (b.1) and (c.1). We note that the reason for this ringing artifact is neither the Gibbs phenomenon nor the decoupling of two subgrids but the fact that the least-squares derivative cannot represent arbitrary deformations of the voxels as discussed in Sec.~\ref{sec:least-squares}.
This discussion, and the outcomes from (b) to (e), indicate that a symmetric and ringing free stress field cannot be obtained by a method with a single deformation gradient per voxel.

Therefore, we also investigated projection operators evaluated on two evaluation points per voxel.
For the Fourier type projection operator on two evaluation points per voxel, as described in \ref{app:fourier-G-2-points-per-voxel}, we still observe ringing, see panel (f). However, ringing is reduced with respect to the Fourier derivative on a single evaluation point. In panel~(f.2) the continuous line represents the stress values in the lower triangle $(P_{xy}^{(1)})$ and the dashed line represents the stress in the upper triangle $(P_{xy}^{(2)})$ (cf. Fig.~\ref{fig:single-element}c). In difference to the stress fields computed with a single evaluation point per voxel the two evaluation points per voxel slightly break the symmetry of the problem. This can be seen in the non-zero stress along the row of the inhomogeneity shown by the red line in (f.2).

Finally, we find a ringing free stress field for the discrete projection operator obtained from linear finite elements on two elements per voxel, panel~(g). The asymmetry between the two elements of a voxel discussed in the previous paragraph persists for this formulation. The stress field result presented in panel (g) seems to be the most appropriate solution to the problem due to its smooth, ringing free field. This conclusion will be supported by the following three examples. We would like to again emphasize that ringing in this example does not only result from the Gibbs phenomenon, which does not exist for a discrete projection operator (as is evident for the forward differences in panel~(d)). A description of local deformation with too few degrees of freedom as shown (e.g., the least square type projection in (e)) also gives rise to ringing.

\subsection{Two pillars and vacuum}

\begin{figure*}
    \centering
    \includegraphics[width=\textwidth]{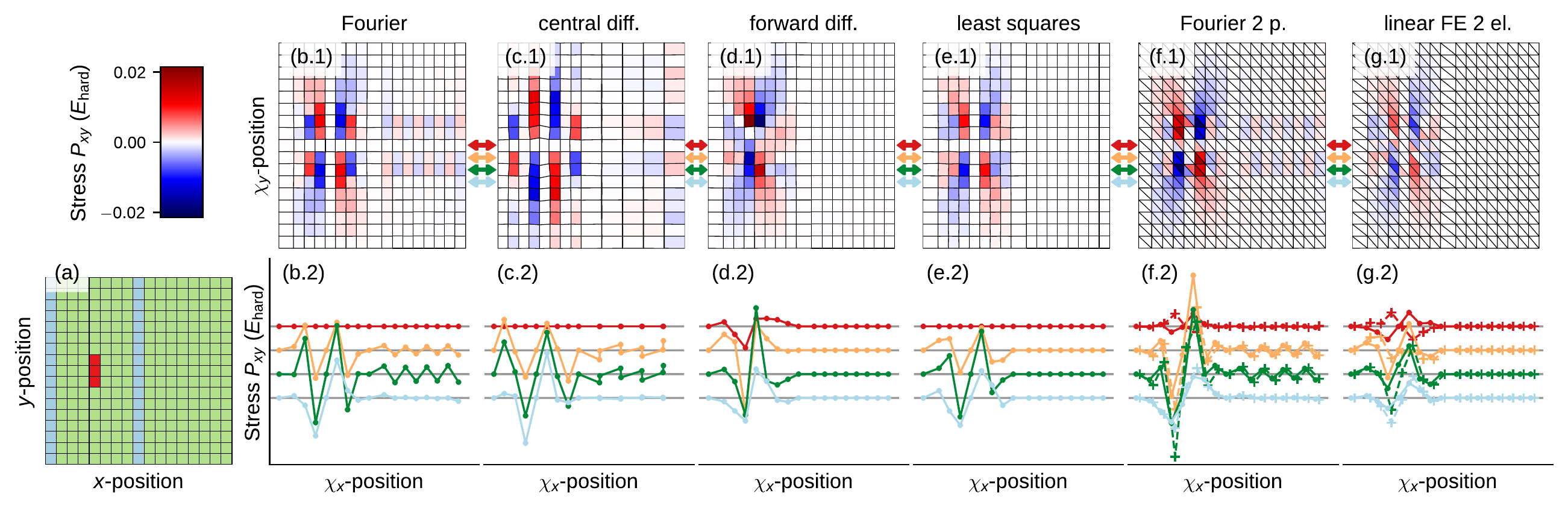}
    \caption{The first Piola-Kirchhoff shear stress $P_{xy}$ of two pillars at 10\% strain in $y$-direction. \textbf{(a)} displays the phase setup. The left pillar (green) has an inhomogeneity (red) of three voxels in its center which are ten times softer than the rest of the pillar. The right pillar (green) is separated by two layers of zero stiffness material (blue) from the left pillar. The first row \textbf{(x.1)} presents the shear stress field for the full grid. The second row \textbf{(x.2)} presents the shear stress as function of $x$ for selected rows of the simulation grid indicated by colored arrows in the first row \textbf{(x.1)}. The rest of the figure is illustrated as described for Fig.~\ref{fig:single_voxel_inhomo}.
    }
    \label{fig:ex_two_pillars}
\end{figure*}

We use a setup consisting of two pillars, as shown in Fig.~\ref{fig:ex_two_pillars}.a, to qualitatively investigate an infinite material contrast and the ability to represent a free surface. This would allow breaking the periodic boundary conditions which are intrinsic to FFT-based methods. We choose a simulation domain of $17 \times 17$ voxels and subdivide it into two pillars (in green) with Youngs modulus $E_{\text{hard}}$ and Poisson number $\nu=0.33$. The two pillars are separated by a layer consisting of single voxel of a material with zero stiffness (in light blue), i.e. $E_{\text{vac}}=0$. One can think of this material as ``air'' or ``vacuum''. At the surface of the pillars we thus have an infinite material contrast. The left pillar, of width $7$ voxels, has additionally an inhomogeneity (in red) of three voxels in its center. The soft inhomogeneity has a Youngs modulus of $E_{\text{soft}}=E_{\text{hard}}/10$ and the same Poisson number of $\nu=0.33$ as the pillar. The inhomogeneity was introduced to generate a non-homogeneous strain field with the ringing artifacts documented in the previous section. The setup is strained by 10\% in the $y$-direction and held at constant size in $x$-direction.
Simulations use the same finite strain model as those of Sec.~\ref{sec:single_pixel_inhomo}.

Figure~\ref{fig:ex_two_pillars} is organized in the same manner as Fig.~\ref{fig:single_voxel_inhomo}. In the first row, panels (b.1) to (g.1), we show the color coded stress field $P_{xy}$. In the second row, panels (b.2) to (g.2), we present the stress field as function of $\chi_x$ at fixed positions $\chi_y$, starting from the center of the inhomogeneity (red line) and going row by row down in orange, green and light blue. The gray vertical lines represent zero stress for each curve. The continuous and dashed lines in panels (f.2) and (g.2) represent the stress of the lower and upper triangular element of a voxel, respectively. The columns are ordered as in Fig. \ref{fig:single_voxel_inhomo}: (b) Fourier type projection, (c) central differences, (d) forward differences, (e) least squares projection, (f) Fourier-type projection on two evaluation points and (g) linear finite element projection on two elements.

For the Fourier-type projection operator in panel (b), we observe ringing artifacts in the left pillar originating from the inhomogeneity. These artifacts are transmitted through the vacuum region to the right pillar. The shear component of the stress should be zero in the right pillar but shows clear ringing artifacts. The ``vacuum'' region of zero stiffness is therefor not able to decouple the two pillars. At the symmetry axis in $x$- and $y$-direction of the inhomogeneity we can again observe a region of zero stress of single voxel thickness.

Central differences, panel (c), lead to a strong decoupling of the grid into two sub grids. The vacuum cannot decouple the strain field in the pillars because the stencil has a range of three voxels. This leads to strong oscillations in $x$-direction in the two subgrids. In the right pillar, that has a width of an even number of grid points in the $x$-direction, the decoupling results in almost zero width of the voxels of one sub grid. (Only four voxels are clearly visible in panel (c.1).) Remarkably the left pillar has no ringing in $y$-direction in the columns of non zero stress. In this example, it becomes very clear that the central differences are sensitive to the setup and number of grid points. For slightly different widths of the pillars or a different mesh grid one can observe these artifacts also in the other pillar (not shown).

Panels (d) and (e) (forward differences and least-squares), show the same behavior as discussed for a single voxel inhomogeneity in the previous Sec.~\ref{sec:single_pixel_inhomo} with an asymmetric, but non-ringing response for the forward differences and a symmetric and ringing response for the least squares scheme. However, the right pillar shows zero shear stress, indicating that for these projection operators the vacuum layer is able to decouple the two pillars. The vacuum layer grows in $x$-direction to absorb the shrinkage of the two pillars (since $\nu=0.33$) while the average strain in $x$-direction is zero. At the surfaces of the left pillar the shear stress field goes to zero as one would expect it for a surface. These discrete derivative schemes decouple the two pillars because their stencil extend only between neighboring nodes.

The Fourier type projection on two evaluation points in panel (f) leads to similar results in the left pillar as for the single voxel inhomogeneity. A ringing artifact from the left pillar is observed in the right pillar which indicates a coupling between the two pillars through the vacuum region. As in the first example, Fig.~\ref{fig:single_voxel_inhomo}f, the two evaluation points lead to a slight symmetry breaking resulting in non-zero stress at the symmetry axis through the center of the inhomogeneity in the left pillar, i.e. the red line in panel (f.2).

For the linear finite element projection shown in panel (g), we again observe artifact-free results. The left pillar shows a smooth, ringing free and symmetric (besides the previous discussed slight asymmetry between upper and lower triangle) stress field originating from the inhomogeneity. The vacuum regions grow in $x$-direction to absorb the shrinking of the pillars in that direction. In the right pillar we find no artifacts from the stress field of the left pillar; thus the pillars are fully decoupled.

This example shows that it is possible to simulate infinite material contrast with vacuum as the soft phase. A single layer consisting of material of zero stiffness can decouple different regions in the RVE and thus break the intrinsic periodic boundary conditions of the FFT-based method in one direction. As expected from the theoretical considerations in Sec.~\ref{chap:methods} the linear finite elements projection operator has the best performance and results in a stress field that is artifact free and qualitatively correct. To further investigate the methods developed here, we continue with a quantitative analysis of an Eshelby inhomogeneity.

\subsection{Eshelby inhomogeneity}

\begin{figure*}
    \centering
    \includegraphics[width=\textwidth]{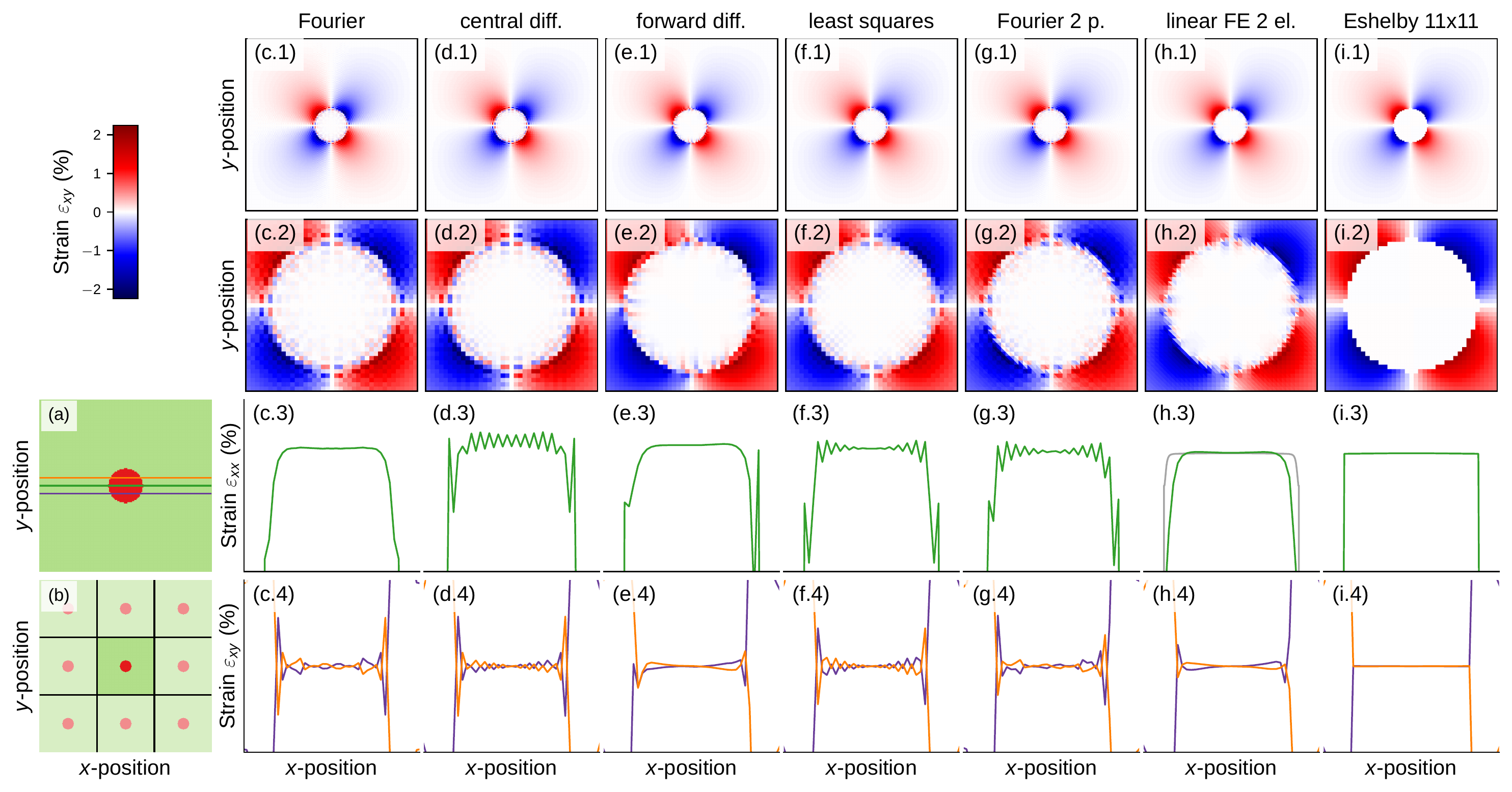}
    \caption{Shear strain of an Eshelby inhomogeneity at biaxial applied strain $\varepsilon_{xx} = \varepsilon_{yy} = 1\%$ at the boundaries on a $151\times151$ grid. \textbf{(a)} shows the phase setup of a soft inhomogeneity (red) in a hard, $E_{\text{hard}}=10E_{\text{soft}}$, matrix (green). The colored lines indicate the rows for which the strain is shown in the third and second row. \textbf{(b)} illustrates the correction of the analytical Eshelby inhomogeneity for periodic boundary conditions. The strain field of the central inhomogeneity is corrected by adding up the strain field from surrounding inhomogeneities as illustrated here in a $3 \times 3$ matrix. The analytical solution presented in column \textbf{(i)} was done for a $11\times 11$ matrix, i.e., including five periodic images in both the positive and negative horizontal and vertical direction. In the first row \textbf{(x.1)} we present the color coded shear strain $\varepsilon_{xy}$ on the whole grid. The second row \textbf{(x.2)} is a zoom on the shear strain field at the inhomogeneity with the same color coding. The third and fourth row show the strain as function of $x$ along selected rows of the zoomed region in row two. The third row \textbf{(x.3)} shows the normal strain in $x$-direction along the green row indicated in \textbf{(a)}. Row four \textbf{(x.4)} presents the shear strain along the middle row of the upper half cylinder (orange) and the middle row through the lower half cylinder (purple) as indicated in \textbf{(a)}. The columns present the results gained with the indicated projection operators: \textbf{(c.i)} Fourier, \textbf{(d.i)} central differences, \textbf{(e.i)} forward differences, \textbf{(f.i)} least squares, \textbf{(g.i)} Fourier on two evaluation points per voxel, \textbf{(h.i)} linear finite elements on two elements per voxel and the column \textbf{(i)} presents the analytical Eshelby solution. In the columns \textbf{(g)} and \textbf{(h)} the third and fourth row show the strain values in the lower triangle, triangle one in Fig.~\ref{fig:single-element}c. In panel \textbf{(h.3)} the additional gray line shows the convergence of the numerical simulation towards the analytical result for an eleven times finer grid with $1661\times1661$ grid points.}
    \label{fig:ex_eshelby_inhomogeneity}
\end{figure*}

The Eshelby inhomogeneity is similar to the first example of a minimal inhomogeneity consisting of a single voxel. The Eshelby inhomogeneity is an ellipsoidal body inside an infinite elastic medium where the elastic medium differs in its material properties from the ellipsoidal body. The analytical solution to the Eshelby problem is well known \cite{Eshelby1957,Eshelby1959,book-Mura1982,Meng2012}. We choose the specific (cylindrical) geometry shown in Fig.~\ref{fig:ex_eshelby_inhomogeneity}a with a hard matrix (light green) of Youngs modulus $E_{\text{hard}}$ and a soft inhomogeneity (red) of Young's modulus $E_{\text{soft}}=E_{\text{hard}}/10$ and zero eigenstrain. The numerical calculations employ a fine mesh of $151\times151$ voxels to properly resolve the cylindrical inhomogeneity with half axes of $10\%$ of the domain edge lengths. The inhomogeneity is placed in the center of the domain and centered on a voxel to retain a symmetric discretized area. For the numerical calculations we use the small-strain formulation (see \ref{sec:small-strain}) since the analytical Eshelby expressions are also obtained in this limit.

The results of these calculations are summarized in Figure~\ref{fig:ex_eshelby_inhomogeneity}. Rows (c.1) to (i.1) show the full solution of the shear strain $\varepsilon_{xy}$ on the $151\times 151$ grid. The next row, panels (c.2) to (i.2), show a zoom of the region containing just the inhomogeneity.
The three colored lines at the center (dark green), upper half (orange) and lower half (purple) of the inhomogeneity in panel (a) indicate the location of the strain components that are shown in the third row, panels (c.3) to (i.3), for the normal strain $\varepsilon_{xx}$ in $x$-direction at the center and in the fourth row, panels (c.4) to (i.4), for the shear strain $\varepsilon_{xy}$ for the upper and lower half of the inhomogeneity. Note that this data is shown only over the zoomed region, not the full calculation. The columns (c) to (h) again represent results obtained for different projection operators as indicated in each column: (c) is the Fourier-type projection, (d) central differences, (e) forward differences, (f) the least square type projection, (g) the Fourier type projection on two evaluation points per voxel and (h) the linear finite elements projection from sec. \ref{sec:linear_finite_elements}.

The additional column (i) represent the analytical results of the Eshelby inhomogeneity. The analytical result is obtained for an inhomogeneity in an infinite media at $1\%$ biaxial strain ($\varepsilon_{xx}=\varepsilon_{yy}=0.01,\ \varepsilon_{xy}=0$) and is corrected for periodic boundary conditions by summing up the influence of $11\times 11$ non-interacting periodic Eshelby inhomogeneities (see panel (b) of Fig.~\ref{fig:ex_eshelby_inhomogeneity}). This correction scheme for periodic boundary conditions converges quickly with the number of images. The average strain on the central inhomogeneity is used as the boundary conditions for the periodic numerical computations.

The Fourier-type projection operator gives qualitatively correct results. Panels (c.3) and (c.4) show a quantitative view of two components of the strain tensor that both clearly show oscillations within the inhomogeneity. As expected, we observe strong ringing artifacts which lead to a deviation from the analytical result especially within the inhomogeneity. Note that the normal strain along the center line vanishes due to symmetry reasons and hence does not show ringing.

For the central differences scheme in panels (d.i), we observe a symmetric checkerboard pattern of decreasing amplitude when approaching the center of the inhomogeneity. At the boundary of the inhomogeneity, there is a double ring-like pattern in the strain field originating from the local decoupling into two subgrids by this scheme.

The forward differences scheme in column (e) produces ringing-free fields but with the drawback of the already discussed asymmetry, best shown in panel (e.3). However the asymmetry of the field is partly suppressed by working in the small strain limit where $\varepsilon_{xy}=\varepsilon_{yx}$. The asymmetry of the strain field is also noticeable in panel (e.4).

Column (f) shows the results produced by the least squares projection operator. As for the two previous examples we observe a checkerboard like pattern of decreasing intensity when approaching the center of the inhomogeneity (see panel (f.2)). 

The Fourier-type projection on two evaluation points per voxel (panel (g)) produces a strain field similar to the standard case with a single point per voxel (panel (c)), however the ringing artifacts appear distributed more homogeneously over the entire inhomogeneity. Panels (g.3) and (g.4) show for clarity only the strain field for the lower triangular element (element one in Figure~\ref{fig:single-element}c). The ringing artifacts in (g.3) and (g.4) are less symmetric than the one of panels (c.3) and (c.4), originating from the symmetry breaking by the triangular mesh.

For the linear finite element projection on two elements shown in column (h), we find ringing free fields and the sharpest drop of the strain at the boundary of the inhomogeneity. Panel (h.3) shows the smoothest curves that are close to the analytical solution presented in (i.3). The normal strain shows a small variation across the inhomogeneity while the analytical solution (panel (i.3)) is almost constant. We find similar variation in the shear strain shown in panel (h.4). The gray line in panel (h.3) demonstrates exemplary the convergence of the numeric simulation towards the analytical result for a eleven times finer grid with $1661 \times 1661$ grid points. Yet, the curves shown in panels (h.3) and (h.4) are closest to the analytical result. 

In summary, we find reasonable agreement with the analytical solution for all projection operators. Linear finite elements gives the smoothest curves and results closest to the analytical findings. As in the previous cases, only the forward differences projection and the finite elements projection eliminate the ringing artifact. For a finer grid we observe a convergence towards the analytical result. It is worth noting that the similar behaviour of the finite elements and forward differences projections is easily explained by the fact that the latter corresponds to the finite element projection where only the lower left triangles are considered.

\subsection{Damage problem}
As a drawback of FFT-based solution methods, ringing artifacts can have a drastic effect on the solution of a homogenization problem. Damage mechanics  problems are especially vulnerable to fluctuations in the stress field caused by ringing artifacts, since localization is one of the most important characteristics of such problems.  A reliable, fast, and ringing-free  homogenization method is therefore essential to address damage mechanics problems.

\begin{figure*}
    \centering
    \includegraphics[width=0.8\textwidth]{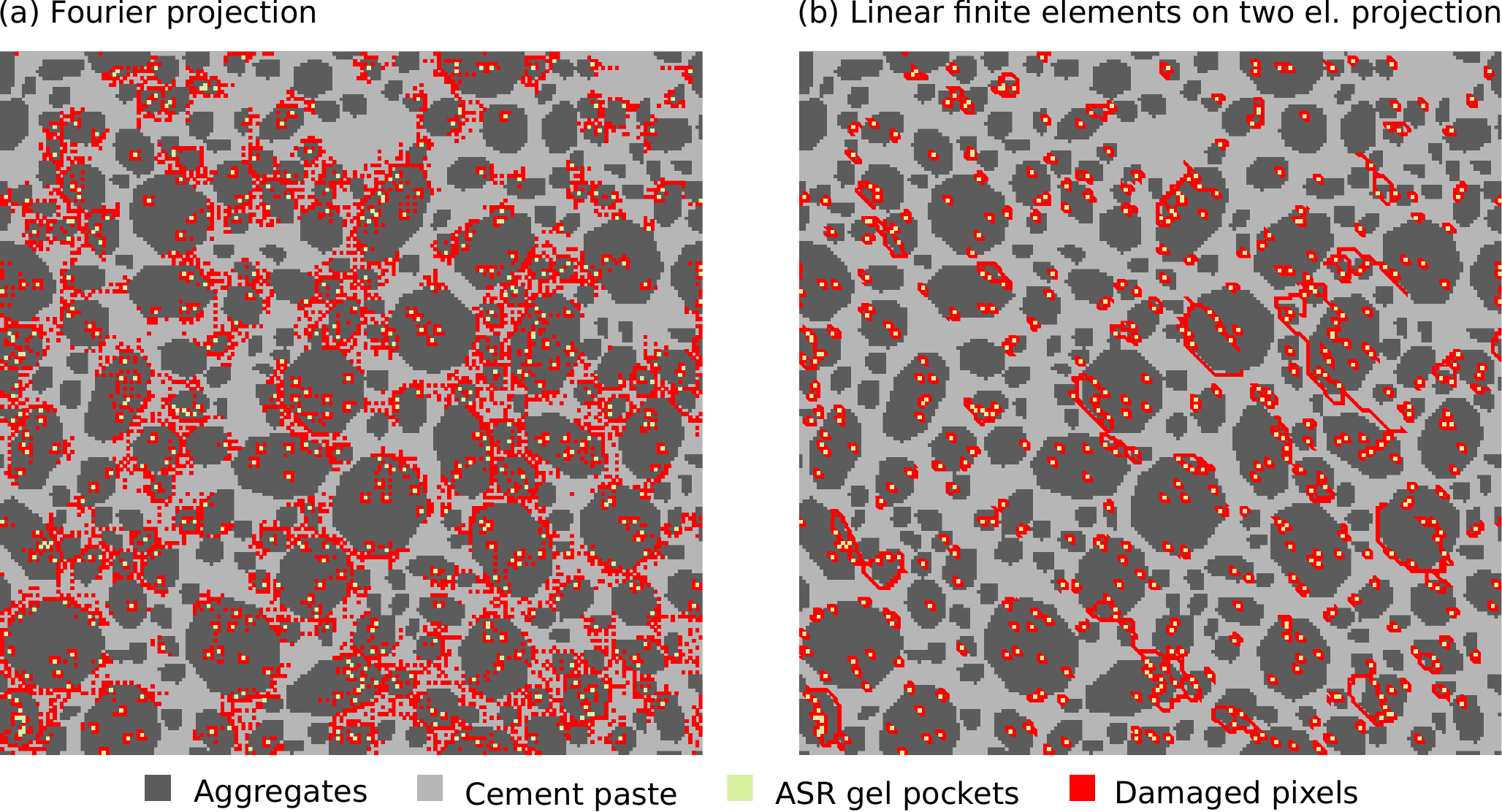}
    \caption{Damage fields in concrete micro-structure RVEs using (a) Fourier projection solver (b) FEM projection solver.}
    \label{fig:damage}
\end{figure*}

In order to illustrate the error introduced by the ringing artifact into a damage problem, we solve a  two-dimensional problem representing a concrete microstructure using an Alkali-Silica reaction (ASR) damage model. In this model, the expansion of gel pockets inside aggregates damages the microstructure.

In the modeled concrete microstructure, three phases are considered — a soft matrix with damage (Cement paste), hard inclusions with damage (Aggregates), and gels whose expansion is modeled by a growing spherical eigenstrain. Polygon-shaped aggregates are placed in the RVE using the level set approach (LSA) algorithm \cite{Sonon2012, Sethian1996}, with an aggregate size distribution chosen to match  sieve sizes from real concrete structures. The aggregate size distribution is truncated on the lower end in order to keep the shape of the aggregates physically sound considering the descritisation grid. ASR gel pockets are placed randomly within the aggregate to fill $2\%$ of the cell surface. Concrete paste and aggregate are represented by a linear damage model as their constitutive law \cite{Ramos2018, Davila2009, Najafgholipour2017}. Due to the brittleness of concrete, the damage part of the bi-linear damage laws is taken steeper than its elastic part. 

In the damage phase, the damage surface threshold is defined by the magnitude of strain measured by the $L_2$-norm. As long as the damage material's strain is below a determined $\epsilon_u$, it behaves as a linear elastic isotropic material with Young modulus of $E_0$ , and afterwards its stress decreases with equivalent stiffness of $-\alpha E_0$ until the stress becomes zero. From that point on, the material does not carry any stress (complete failure). An eigenstrain with the final amplitude of $20 \epsilon_u$ was applied on gel voxels, placed as explained beforehand, in $1000$ consecutive steps in the carried out simulation.
The damage field caused by this loading scenario is depicted in Fig.~\ref{fig:damage} employing Fourier and FEM projection solvers. As observed in Fig.~\ref{fig:damage} the damage pattern evolved in the Fourier projection solver solution is checker-boarded and therefore non-physical. While as demonstrated in Fig.~\ref{fig:damage}b, after damage initiation around ASR gel pockets, micro-cracks formed during the damage process coalesce to form cracks with lengths in the range of 0.2 times the RVE length. Comparison of Fig.~\ref{fig:damage}a, b suggests that, in contrast with Fourier projections solvers, ringing-free spectral FEM projection solvers are capable to simulate mechanics damage RVE problems.

\section{Summary \& Conclusion}

We have extended the compatibility projection of Vond\v{r}ejc, Zeman, de Geus and coworkers~\cite{Vondrejc2014,Zeman2017,deGeus2017} to basis sets with local support. The nonlocal support of the Fourier basis gives rise to ringing artifacts and makes it impossible to include regions with zero stiffness or vacuum into a calculation. We show that using linear finite elements as a basis set eliminates all of the problems of the Fourier basis but retains the advantages of the compatibility projection, such as the rapid convergence rate. In particular, this formulation allows modeling of regions with zero stiffness, which opens the method to the study of free surfaces or metamaterials. We note that the ideas behind compatibility projection can be exploited to construct preconditioners for standard finite-element formulations beyond the linear elements presented here~\cite{ladecky_well-conditioned_nodate}.

Among the projection operators examined in this paper, all but the least squares operator guarantee a compatible strain field. It is useful to think of those compatible projection operators in two categories; discretisations with local support and discretisations with global support, where local support refers to non-overlapping projection stencils. We observe that only two local support projections, the finite element projection and the forward differences projection eliminate ringing, unlike the Fourier and central differences operators (global support and non-local support including neighbouring voxels). This observation allows a well known property of finite element calculations, where discretisation errors are to be expected to be significant at the length scale of the element size (= support size) by St. Venant's principle. This observation allows to interpret the ringing artifact as a discretization error to be expected at the length scale of the support.

\section*{Acknowledgments}

We acknowledge funding by the European Research Council (StG-757343~(LP)), the Carl Zeiss Foundation (Research cluster "Interactive and Programmable Materials - IPROM"~(LP)), the Deutsche Forschungsgemeinschaft (EXC 2193/1 - 390951807~(LP)), the Swiss National Science
Foundation (Ambizione grant 174105~(TJ)), the Czech Science Foundation (projects No.~20-14736S~(ML) and 19-26143X~(JZ)), and the European Regional Development Fund (Centre of Advanced Applied Sciences – CAAS, CZ.02.1.01/0.0/0.0/16\_019/0000778 (ML, IP)).

\appendix

\section{Discrete Fourier transformation}
\label{app:discrete-fourier-trafo}

The 2D discrete Fourier transformation is defined as follows through out the paper. The generalization to 1D or 3D is straight forward. We divide the simulation domain of edge lengths $L_x, L_y$ into $N_x, N_y$ voxels of equal edge lengths $\Delta_i = L_i/N_i$ in each spatial direction $i\in\{x,y\}$. The lower left corner of voxel $n_x=I, n_y=J$ is then given by
\begin{equation}
    r_{IJ,i} = \frac{L_i}{N_i} n_i, \quad n_i =
        0, 1, \dots, N_i - 1.
\end{equation}
The corresponding wave vectors $\{\vec{k}_{IJ}\}$ with $m_x=I, m_y=J$ are
\begin{align}
    \label{eq:Fourier_wave_vectors}
    k_{IJ,i} =&  \frac{2\pi}{L_i} m_i,\\
    \quad m_i =& 
    \begin{cases}
        -\frac{N_i}{2},\dots, 0, 1, \dots, \frac{N_i}{2} - 1 \  , \quad N_i \text{ even}, \\
        -\frac{N_i-1}{2}, \dots, 0, 1, \dots, \frac{N_i-1}{2} \ , \quad N_i \text{ odd}.
    \end{cases} \nonumber
\end{align}
For the discrete Fourier transform of the function $\vec{f}(\vec{r})$ we use
\begin{equation}
    \label{eq:discrete-fourier-trafo}
    \mathcal{F}\left(\vec{f}(\vec{r})\right)(\vec{k}) = \hat{\vec{f}}(\vec{k}) = \sum_{\vec{r}} \vec{f}(\vec{r}) \exp \left(- \ii \vec{k} \cdot \vec{r} \right)
\end{equation}
with the corresponding inverse transformation
\begin{equation}
    \label{eq:inverse-discrete-fourier-trafo}
    \mathcal{F}^{-1}\left(\hat{\vec{f}}(\vec{k})\right)(\vec{r}) = \vec{f}(\vec{r}) = \frac{1}{N} \sum_{\vec{k}} \hat{\vec{f}}(\vec{k}) \exp \left(\ii \vec{k} \cdot \vec{r} \right).
\end{equation}

\section{Small-strain projection}
\label{sec:small-strain}

It is often useful to carry out calculations in the small-strain limit. Small strains are special because the strain tensors (that replaces the deformation gradient) has to remain symmetric. For a formulation that involves multiple elements per voxel, we also require multiple symmetric strain tensors per voxel. While in the finite strain formulation, these were absorbed in our derivative indices $\alpha$, $\beta$ etc., we need to introduce a specific element index for the small-strain case and can no longer distinguish between derivatives and Cartesian coordinates because of the symmetry in between derivatives and coordinates. Additionally to the previous introduced indices we will therefore use capital Greek letters ($\Theta,\, \Lambda,\, \Xi,\,\dots$) to denote elements. Components of the strain tensor will be denoted by small Latin indices.

We now introduce the strain tensor $\t{e}$ in lieu of the deformation gradient, Eq.~\eqref{eq:def-grad}. The strain tensor for element $\Theta$ is
\begin{equation}
    \t{e}_\Theta = \frac{1}{2}\left[\nabla_\Theta\otimes\vec{u} + \left(\nabla_\Theta\otimes\vec{u}\right)^T\right],
    \label{eq:small-strain}
\end{equation}
where $\nabla_\Theta$ is the (potentially discrete) derivative operator for element $\Theta$ and $\vec{u}(\vec{r})=\vec{\chi}(\vec{r})-\vec{r}$ are the displacements from the undeformed positions $\vec{r}$.
For a given $\t{e}$, we minimize the residual $\mathcal{R}=\sum_\Theta \sum_{\vec{k}} \t{R}^*_\Theta(\vec{k}):\t{R}_\Theta(\vec{k})$ with
\begin{equation}
    \t{R}_\Theta(\vec{k}) = \frac{1}{2}\left(\hat{\vec{\mathcal{D}}}_\Theta(\vec{k}) \otimes \hat{\vec{u}}(\vec{k}) + \hat{\vec{u}}(\vec{k}) \otimes \hat{\vec{\mathcal{D}}}_\Theta(\vec{k}) \right) - \hat{\t{e}}_\Theta(\vec{k})
\end{equation}
with respect to $\vec{u}^*$,
where we have transformed into the Fourier-space and introduced the gradient operator $\hat{\vec{\mathcal{D}}}_\Theta(\vec{k})$.
The full residual is given by
\begin{equation}
\begin{split}
    2\mathcal{R}
    =
    \sum_\Theta \bigg(
    &
    \hat{\vec{\mathcal{D}}}_\Theta^*\cdot \hat{\vec{\mathcal{D}}}_\Theta
    \,
    \hat{\vec{u}}^*\cdot \hat{\vec{u}}
    +
    \hat{\vec{\mathcal{D}}}_\Theta^*\cdot \hat{\vec{u}}
    \,
    \hat{\vec{u}}^*\cdot \hat{\vec{\mathcal{D}}}_\Theta
    \\
    &
    - \hat{\vec{\mathcal{D}}}_\Theta^* \cdot \t{e}_\Theta\cdot \hat{\vec{u}}^*
    - \hat{\vec{u}}^* \cdot \t{e}_\Theta\cdot \hat{\vec{\mathcal{D}}}_\Theta^* \\
    &
    - \hat{\vec{\mathcal{D}}}_\Theta \cdot \t{e}_\Theta^* \cdot \hat{\vec{u}}
    - \hat{\vec{u}} \cdot \t{e}_\Theta^* \cdot \hat{\vec{\mathcal{D}}}_\Theta + 2\t{e}^*_\Theta:\t{e}_\Theta
    \bigg).
\end{split}
\end{equation}
Minimization yields
\begin{equation}
    \left(\t{1} + \sum_\Theta \hat{\t{g}}_{\Theta\Theta}\right)\cdot\hat{\vec{u}}
    =
    \sum_\Theta\left(
    \hat{\vec{\mathcal{D}}}_\Theta^{-1}\cdot\hat{\t{e}}_\Theta
    +
    \hat{\t{e}}_\Theta\cdot\hat{\vec{\mathcal{D}}}_\Theta^{-1}
    \right)
    \label{eq:small-strain-eq2}
\end{equation}
with
\begin{equation}
    \hat{\vec{\mathcal{D}}}_\Theta^{-1} = \frac{\hat{\vec{\mathcal{D}}}_\Theta^*}{\sum_\Lambda \hat{\vec{\mathcal{D}}}_\Lambda \cdot \hat{\vec{\mathcal{D}}}_\Lambda^*}
    \quad\text{and}\quad
    \hat{\t{g}}_{\Theta\Lambda} = \hat{\vec{\mathcal{D}}}_\Theta\otimes\hat{\vec{\mathcal{D}}}_\Lambda^{-1}.
\end{equation}
Equation~\eqref{eq:small-strain-eq2} can be formally solved to give
\begin{equation}
    \hat{\vec{u}} = \hat{\t{h}}\cdot\sum_\Theta\left(\hat{\vec{\mathcal{D}}}_\Theta^{-1}\cdot\hat{\t{e}}_\Theta+\hat{\t{e}}_\Theta\cdot\hat{\vec{\mathcal{D}}}_\Theta^{-1}\right).
    \label{eq:small-strain-chi}
\end{equation}
or
\begin{equation}
    \hat{u}_i = \sum_\Theta\left(\hat{D}^{-1}_{\Theta,l} \hat{h}_{im} + \hat{h}_{il}\hat{D}^{-1}_{\Theta,m}\right)\hat{e}_{\Theta,lm}.
\end{equation}
with
\begin{equation}
    \hat{\t{h}} = \left(\t{1} + \sum_\Theta \hat{\t{g}}_{\Theta\Theta}\right)^{-1}.
\end{equation}
To arrive at the projection operator, we now need to insert Eq.~\eqref{eq:small-strain-chi} into Eq.~\eqref{eq:small-strain}. This yields the small-strain projection operator
\begin{equation}
    \hat{G}_{\Theta\Lambda,ijlm}
    =
    \frac{1}{2}\bigg(
    \hat{g}_{\Theta\Lambda,il} \hat{h}_{jm}
    +
    \hat{g}_{\Theta\Lambda,im} \hat{h}_{jl}
    +
    \hat{g}_{\Theta\Lambda,jl} \hat{h}_{im}
    +
    \hat{g}_{\Theta\Lambda,jm} \hat{h}_{il}
    \bigg),
\end{equation}
where the projected strains are given by
\begin{equation}
    \hat{e}_{\Theta,ij} = \sum_\Lambda \hat{G}_{\Theta\Lambda,ijlm} \hat{e}_{\Lambda,ml}.
\end{equation}
By combining the pairs of indices $\alpha=\Theta,j$ and $\beta=\Lambda,l$, we can write this in the same form as the large strain projection, $\hat{e}_{i\alpha}=\hat{G}_{i\alpha\beta j} \hat{e}_{j\beta}$.
Note that for a single element $\Theta=1$, we can write down the expression for $\hat{\t{h}}$ analytically,
\begin{equation}
    \hat{\t{h}} = \t{1} - \frac{1}{2} \hat{\t{g}}_{11}.
\end{equation}
Using this expression and the Fourier derivative for $\hat{\vec{\mathcal{D}}}(\vec{k})$ gives the small-strain projection operator of Ref.~\cite[Section~6]{Milton1988}.

\section{Fourier-type projection operator on two evaluation points per voxel}
\label{app:fourier-G-2-points-per-voxel}

The Fourier-type derivative can be extended to several gradient evaluation points per voxel. This is useful to investigate the influence of Gibbs ringing separately from the effect of missing degrees of freedom. The standard Fourier-type derivative is evaluated at the grid points of the Fourier grid, which are the centers of the voxels. Hence, for a triangular mesh we additionally evaluate the Fourier derivative at the geometrical center of each triangle. For a two dimensional rectangular cell of edge lengths $\Delta_i$, as shown in Fig.~\ref{fig:single-element}c, that means applying a shift of $\pm ( \Delta_1 , \Delta_2 )/6$ from the center. The Fourier derivative operator $\hat{\mathcal{D}}_\alpha(\vec{k}) = \ii k_\alpha$ acquires a phase to yield
\begin{align}
\begin{split}
    \hat{\mathcal{D}}_{1,i}(\vec{k}) = \ii k_i \exp \left( \frac{-\ii}{6} \sum_i k_i \Delta_i \right), \\
    \hat{\mathcal{D}}_{2,i}(\vec{k}) = \ii k_i \exp \left( \frac{+\ii}{6} \sum_i k_i \Delta_i \right),
\end{split}
\end{align}
where we used explicit the indices $1$ and $2$ to denote the derivative operator in the center of the lower and upper triangle as shown in Fig.~\ref{fig:single-element}c.

\end{document}